\documentclass[11pt]{article}

\usepackage{epsfig,graphics}
\usepackage{amssymb,amsmath,amsthm,bm}

\newtheorem{theo}{Theorem}
\newtheorem{lem}{Lemma}

\newtheorem{cor}{Corollary}
\newtheorem{mrem}{Remark}

\def\R{{\mathbb R}}

\def\ga{\alpha}
\def\gga{\gamma}

\def\go{\omega}

\def\gs{\sigma}
\def\gl{\lambda}
\def\wt{\widetilde}

\def\gp{{\prime}}

\def\ep{\epsilon}
\def\vep{\varepsilon}

\def\gt{\triangle}

\def\b0{{\bf 0}}
\def\1{{\bf 1}}

\def\bd{{\bf d}}

\def\cS{\mathcal {S}}

\def\rem{{\rm rem}}
\def\Rem{{\rm Rem}}

\def\err{{\rm res}}
\def\Err{{\rm Res}}
\def\bd{{\rm bd}}
\def\Bd{{\rm Bd}}

\def\pd{{\partial}}

\def\wh{\widehat}

\def\wt{\widetilde}

\begin{document}

\title{Analysis of Adaptive Short-time Fourier Transform-based Synchrosqueezing Transform\\
{\small November 28, 2018}
}

\author{Haiyan Cai${}^{1}$,  Qingtang Jiang${}^1$, Lin Li${}^1$, and Bruce W. Suter${}^2$
}

\date{}

\maketitle



\bigskip 

{\small 1. Dept. of Math \& CS, University of Missouri-St. Louis, St. Louis,  MO 63121, USA

}

\medskip 

{\small 2. The Air Force Research Laboratory, AFRL/RITB, Rome, NY 13441, USA.  

}

\begin{abstract}
Recently the study of modeling a non-stationary signal as
a superposition of 
amplitude and frequency-modulated Fourier-like oscillatory modes has been a very active research area. The synchrosqueezing transform (SST) is a powerful method 
for instantaneous frequency estimation and component separation of non-stationary multicomponent signals.
 The short-time Fourier transform-based SST (FSST for short) reassigns the frequency variable to sharpen the time-frequency representation and to separate the components of a multicomponent non-stationary signal. Very recently the FSST with a time-varying parameter, called the adaptive FSST, was introduced. 
The simulation experiments 
show that the adaptive FSST is very promising in instantaneous frequency estimation of the component of 
a multicomponent signal, and in accurate component recovery.  
However the theoretical analysis of the adaptive FSST has not been carried out. In this paper, 
we study the theoretical analysis of the adaptive FSST and obtain the error bounds for 
the instantaneous frequency estimation and component recovery with the adaptive FSST and the 2nd-order adaptive FSST.

\end{abstract}

\section{Introduction}


To model a non-stationary signal  $x(t)$ as
\begin{equation}
\label{AHM}
x(t)=\sum_{k=1}^K x_k(t), \qquad x_k(t)=A_k(t) e^{i 2\pi \phi_k(t)},  
\end{equation}
with $A_k(t), \phi_k'(t)>0$ 
is important to extract information, such as the underlying dynamics,  hidden in $x(t)$.

The empirical mode decomposition (EMD) \cite{Huang98} is a widely used method for the representation 
of \eqref{AHM}. Recently  
the continuous wavelet transform (CWT)-based  synchrosqueezed transform (WSST) developed in  \cite{Daub_Lu_Wu11} provides a mathematically sound alternative to EMD. 
The short-time Fourier transform (STFT)-based SST (FSST) was proposed in \cite{Thakur_Wu11} and further studied in \cite{Wu_thesis, MOM14}. 
Both WSST and FSST are a special type of the reassignment method \cite{A_Flandrin_reassignment95}. 
WSST reassigns the scale varilable of the CWT  to the frequency variable and  FSST  reassigns the frequency variable, both aiming to sharpen the time-frequency representation and to separate the components of a multicomponent non-stationary signal.
SST was proved to be robust to noise and small perturbations \cite{Thakur_etal_Wu13, Iatsenko15, Meignen17}. However SST does not provide sharp representations for signals with significantly frequency changes. In this regard, 
the  2nd-order FSST and  the 2nd-order WSST 
were introduced in \cite{MOM15} and \cite{OM17} respectively, 
and theoretical analysis of the 2nd-order FSST was carried out in \cite{BMO18}. 
The 2nd-order SST improves the concentration of the time-frequency representation.
Other SST related methods include the generalized WSST \cite{Li_Liang12},
a hybrid EMT-SST computational scheme \cite{Chui_Walt15},
the synchrosqueezed wave packet transform \cite{Yang15},
WSST with vanishing moment wavelets \cite{Chui_Lin_Wu15},
the multitapered WSST \cite{Daub_Wang_Wu15}, the demodulation-transform based SST  \cite{Wang_etal14, Jiang_Suter17, WCSGTZ18},
higher-order FSST \cite{Pham17}, signal separation operator \cite{Chui_Mhaskar15} and empirical signal separation algorithm \cite{LCJJ17}. The statistical analysis of synchrosqueezed transforms has been studied in \cite{Yang18}.

Most of the  WSST or FSST algorithms available in the literature are based on a continuous (admissible) wavelet or a window function with a  fixed window, which means high time resolution and frequency resolution cannot be obtained simultaneously.  For broadband signals, a narrow window is suitable for the high-frequency parts while  a wide window is suitable for the low-frequency parts.  Recently the authors in \cite{Wu17} introduced  a method to select the time-varying window width for sharp SST representation by minimizing the R${\rm \acute e}$nyi entropy. The window width of the signal separation operator  algorithm in \cite{Chui_Mhaskar15} is also time-varying.  
More recently the authors of \cite{Saito17, Berrian18} study the adaptive FSST with the window function containing time and frequency parameters. Very recently the authors of \cite{LCHJJ18, LCJJ18} proposed  
 the adaptive WSST and adaptive FSST 
 with a time-varying adaptive Gaussian window. They obtain the well-separated condition for multicomponent signals  using the linear frequency modulation to approximate a non-stationary signal during any local time, along with a new definition of bandwidth of Gaussian window. The experimental experiments with synthetic and real data show that the adaptive FSST is very promising in instantaneous frequency estimation of the component of 
a multicomponent signal, and in accurate component recovery.  
However the theoretical analysis of the adaptive FSST has not been carried out. 
The goal of this paper is to study the theoretical analysis of such an FSST. 
We obtain the error bounds for the instantaneous frequency estimation and component recovery with the adaptive FSST and the 2nd-order adaptive FSST. 

The rest of this paper is organized as follows. In Section 2 we  briefly review  FSST, the 2nd-order FSST, the adaptive FSST and the 2nd-order adaptive FSST. In Section 3, we obtain the theoretical analysis of the (1st-order) adaptive FSST. We establish the error bounds for the IF estimation and component recovery. In Section 4, we consider the theoretical analysis of the 2nd-order adaptive FSST. The error bounds for the IF estimation and component recovery for the 2nd-order adaptive FSST are obtained.  The proof two lemmas is presented in the appendix.

\section{Short-time Fourier transform-based synchrosqueezed transform}

In this section we briefly review the short-time Fourier transform (STFT)-based synchrosqueezed transform (FSST) and the adaptive FSST. 
The (modified) STFT 
of $x(t)\in L_2(\R)$ with  a window function $g(t)\in L_2(\R)$
is defined by
\begin{eqnarray*}
\label{def_STFT}
V_x(t, \eta) \hskip -0.6cm&& =\int_\R x(\tau)g(\tau-t)e^{-i2\pi \eta (\tau-t)}d\tau, 
\end{eqnarray*}
where $t$ and $\eta$ are the time variable and the frequency variable respectively.  

The original signal $x(t)$  can be recovered back from its STFT:
\begin{equation*}
x(t)=\frac1{\|g\|_2^2}\int_\R\int_\R
V_x(t, \eta) \overline {g(t-\tau)} e^{-i2\pi \eta (\tau-t)}d\tau  d\eta.
\end{equation*}
If $g(0)\not=0$, then one can show that 
$x(t)$  can also be recovered back from its STFT $V_x(t, \eta) $ with integrals involving only $\eta$: 
\begin{equation}
\label{reconst_STFT}
x(t)=\frac1{g(0)}\int_\R
V_x(t, \eta) d\eta.
\end{equation}
In addition, if  the window function $g(t)\in L_2(\R)$ is real, then for a real-valued $x(t)\in L_2(\R)$, we have
\begin{equation}
\label{reconst_STFT_real}
x(t)= \frac 2{g(0)} {\rm Re} \Big( \int_0^\infty V_x(t, \eta) d\eta\Big).
\end{equation}

Here we remark that if the window function $g(t)$ is in the Schwarz class $\cS$,
then STFT  $V_x(t, \eta)$ of a slowly growing $x(t)$ with $g(t)$ is well defined. Furthermore, the above formulas still hold.

In this paper we assume $g(0)\not=0$. For a signal $x(t)$, its Fourier transform $\wh x(\xi)$ is  defined by
$$
\wh x(\xi)=\int_\R x(t)e^{-i2\pi \xi t} dt.
$$


\subsection{STFT-based synchrosqueezing transform}


The STFT-based synchrosqueezing transform (FSST) was first studied in \cite{Thakur_Wu11}. 
For a signal $x(t)$, at $(t, \eta)$ for which $V_x(t, \eta)\not=0$,
denote
\begin{equation*}
\go_x(t, \eta)={\rm Re}\big(\frac{\partial _t V_x(t, \eta)}{2\pi i V_x(t, \eta)}\big). 
\end{equation*}
The quantity $\go_x(t, \eta)$ is called the ``phase transformation"  \cite{Daub_Lu_Wu11} or 
``instantaneous frequency  information" in \cite{Thakur_Wu11}. 
FSST is to reassign the frequency variable $\eta$ by transforming STFT $V_x(t, \eta)$ of $x(t)$ to a quantity, denoted by $R_{x, \gga}^{\gl}(t, \xi)$, on the time-frequency plane:
\begin{equation}
\label{def_FSST}
R_{x, \gga}^{\gl}(t, \xi)
=\int_{|V_x(t, \eta)|>\gga}  V_x(t, \eta) \frac 1{\gl}h\big(\frac{\xi-\go_x(t, \eta)}\gl\big) d\eta. 
\end{equation}
where throughout this paper $h(t)$ is a compactly supported function with certain smoothness and $\int_\R h(t)dt=1$.  Throughout 
this paper $\int_{|V_x(t, \eta)|>\gga}$ means the integral  $\int_{\{\eta: \; |V_x(t, \eta)|>\gga\}}$ with  
$\eta$ over the set $\{\eta: \; |V_x(t, \eta)|>\gga\}$.


We consider multicomponent signals $x(t)$ given by \eqref{AHM}, 
where $A_k(t), \phi_k(t)$ satisfy 
\begin{eqnarray}
 \label{cond_basic0}&&A_k(t)\in C^1(\R)\cap L_\infty(\R), \phi_k(t)\in C^2(\R), \\
\label{cond_basic}&& A_k(t)>0, \; \inf_{t\in R} \phi^\gp_k(t)>0, \; \sup_{t\in R} \phi^\gp_k(t)<\infty.
\end{eqnarray}
Let $\vep>0$ and $\gt>0$, and let ${\cal B}_{\vep, \gt}$ denote the set of multicomponent signals of \eqref{AHM} satisfying \eqref{cond_basic0}, \eqref{cond_basic}, and  the following conditions: 
\begin{eqnarray}
 \label{cond_phi_2nd_der}
 && |A^\gp_k(t)|\le \vep \phi^\gp_k(t),  \;  |\phi^{\gp\gp}_k(t)|\le \vep\phi^\gp_k(t), \; t\in \R, \; 
 M^{\gp\gp}_k=\sup_{t\in \R} |\phi^{\gp\gp}_k(t)|<\infty, \\
\label{freq_resolution}&&\phi^\gp_k(t)-\phi^\gp_{k-1}(t)\ge 2\gt, \; 2\le k\le K, t\in \R. 
\end{eqnarray}
The condition \eqref{freq_resolution} is called the well-separated condition with resolution $\gt$. 
For the  well-separated condition, \cite{Thakur_Wu11} uses a stronger condition 
 than that in \eqref{freq_resolution}:
 \begin{equation}
 \label{freq_resolution_old} \inf_{t\in \R} \phi^\gp_k(t)- \sup_{t\in \R}\phi^\gp_{k-1}(t)\ge 2\gt, \; 2\le k\le K. 
\end{equation}
The condition \eqref{cond_phi_2nd_der}, which was considered in \cite{Wu_thesis}, means that $A_k(t)$ and IF $\phi_k^\gp(t)$ change slowly compared with $\phi_k(t)$. 
\cite{MOM14} uses another condition for the change of   $A_k(t)$ and IF $\phi_k^\gp(t)$ :
\begin{eqnarray}
\label{cond_phi_2nd_der1}&& |A^\gp_k(t)|\le \vep,  \;  |\phi^{\gp\gp}_k(t)|\le \vep, \; t\in \R.  
\end{eqnarray}
We let $B_{\vep, \gt}$ denote the set of multicomponent signals of \eqref{AHM} satisfying \eqref{cond_basic0}, \eqref{cond_basic}, \eqref{cond_phi_2nd_der1} and well-separated condition \eqref{freq_resolution}. 

Let 
\begin{equation}
\label{def_Zk0}
{\cal Z}_k=\{\eta: \; |\eta-\phi^\gp_k(t)|< \gt\}. 
\end{equation}
Then the well-separated condition \eqref{freq_resolution} implies that ${\cal Z}_k, 1\le k\le K$ are not overlapping.

\bigskip

Denote 
\begin{eqnarray}
&&\label{def_Gamma0} \Gamma_0(t)=K I_1 +\pi  I_2 \sum_{k=1}^K A_k(t), \; 
\wt \Gamma_0(t)=K \wt I_1 +\pi  \wt I_2 \sum_{k=1}^K A_k(t), 
\end{eqnarray}
where 
\begin{eqnarray}
&&\label{def_In}
\label{def_tIn}
I_n=\int_\R  | \tau^n g(\tau)| d\tau, \; 
\wt I_n=\int_\R  | \tau^n g^\gp(\tau)| d\tau, \; n=1, 2, \cdots. 
\end{eqnarray}

{\bf Theorem A.} {\it Let $x(t)\in B_{\vep, \gt}$ and  $g$ be a function in the Schwartz class with supp($\wh g)\subseteq [-\gt, \gt]$. 
Let $\Gamma_0(t), \wt \Gamma_0(t) $ be defined by \eqref{def_Gamma0}. 
Then we have the following. 

  {\rm (a)} 
  Suppose $\wt \vep$ satisfies  $\wt \vep \ge \vep \Gamma_0(t)$. Then for any $\eta$ with 
  $|V_x(t, \eta)|>\wt \vep$,  there exists a unique $k\in \{1, 2, \cdots, K\}$ such that $(t, \eta)\in {\cal Z}_k$.   

{\rm (b)} Suppose $(t, \eta)$ satisfies $|V_x(t, \eta)|>\wt \vep$ and  $(t, \eta)\in {\cal Z}_k$. Then 
\begin{equation}
\label{transformation_est0}
|\go_x(t, \eta)-\phi^\gp_k(t)|<  \frac {\vep}{\wt \vep}\big(\Gamma_0(t) \gt +\frac 1{2\pi} \wt \Gamma_0(t)\big).  
\end{equation}

{\rm (c)} Suppose that $\wt \vep$ satisfies 
$ \big(\Gamma_0(t) \gt +\frac 1{2\pi} \wt \Gamma_0(t)\big) {\vep}/{\wt \vep} \le \gt$. 
Then, for any $k\in \{1, \cdots, K\}$ and any  $\wt\vep_3$ satisfying 
  $\big(\Gamma_0(t) \gt +\frac 1{2\pi} \wt \Gamma_0(t)\big) {\vep}/{\wt \vep}\le \wt \vep_3 \le \gt$, we have 
\begin{equation}
 \label{reconstr_1st0}
 \Big| \lim_{\gl\to 0} \frac1{g(0)}\int_{|\xi-\phi^\gp_k(t)|<\wt \vep_3} R_{x, \wt \vep}^{\gl}(t, \xi)d\xi -x_k(t) \Big|\le \frac {2\gt(\vep \Gamma_0(t)+ \wt \vep)}{|g(0)|}.   
\end{equation}

   {\rm (d)} If   $x(t)\in {\cal B}_{\vep, \gt}$, then the above statements {\rm (a)-(c)} hold with  
 $  \Gamma_0(t)$ and $\wt \Gamma_0(t)$  in \eqref{def_Gamma0} replaced by 
  \begin{eqnarray}
&&\label{def_Gamma0_Wu} \Gamma_0(t)=\sum_{k=1}^K\Big\{ \phi_k^\gp(t) I_1 +\frac 12 M_k^{\gp\gp}I_2 +\pi  A_k(t)(
\phi_k^\gp(t) I_2 +\frac 13 M_k^{\gp\gp}I_3)\Big\}, \\
&&\label{def_Gamma0_der_Wu} \wt \Gamma_0(t)=\sum_{k=1}^K\Big\{ \phi_k^\gp(t) \wt I_1 +\frac 12 M_k^{\gp\gp}\wt I_2 +\pi  A_k(t)(\phi_k^\gp(t) \wt I_2 +\frac 13 M_k^{\gp\gp}\wt I_3)\Big\}. 
\end{eqnarray}
  }

Here we remark that $\wt \vep$ and $\wt \vep_3$ in Theorem A could be a function of $t$.  

If we choose $\wt \vep={\vep}^{1/3}$ and if $\vep$ is small enough such that 
\begin{equation}
\label{cond00} 
\wt \vep\le \min\{\gt, \frac 1{\|\Gamma_0(t)\gt +\frac 1{2\pi}\wt  \Gamma_0(t)\|_\infty}\}, 
\end{equation}
then $\wt \vep (\Gamma_0(t)\gt +\frac 1{2\pi} \wt \Gamma_0(t))\le 1$. Hence, 
$$
\big(\Gamma_0(t) \gt +\frac 1{2\pi} \wt \Gamma_0(t)\big) {\vep}/{\wt \vep}\le \wt \vep \le \gt.
$$ 
Thus, the conditions in Theorem A are satisfied,  and Theorem A (with $\wt \vep_3=\wt \vep$) can be stated in the following theorem.

{\bf Theorem B.}  \cite{Wu_thesis, MOM14} {\it Let $x(t)\in {\cal B}_{\vep, \gt}$ or  $B_{\vep, \gt}$, and $\wt \vep=\vep^{1/3}$. Let $g$ be a function in the Schwartz class with supp($\wh g)\subseteq [-\gt, \gt]$. 
 If $\vep$ is small enough, then  the following statements hold. 

  {\rm (a)} 
 For $(t, \eta)$ satisfying  $|V_x(t, \eta)|>\wt \vep$, there exists a unique $k\in \{1, 2, \cdots, K\}$ such that $(t, \eta)\in {\cal Z}_k$.   

{\rm (b)} Suppose $(t, \eta)$ satisfies $|V_x(t, \eta)|>\wt \vep$ and  $(t, \eta)\in {\cal Z}_k$. Then 
\begin{equation*}
|\go_x(t, \eta)-\phi^\gp_k(t)|<  \wt \vep.  
\end{equation*}

{\rm (c)} 
For any $k\in \{1, \cdots, K\}$, 
\begin{equation*}
 \label{reconstr00}
 \Big| \lim_{\gl\to 0} \frac1{g(0)}\int_{|\xi-\phi^\gp_k(t)|<\wt \vep} R_{x, \wt \vep}^{\gl}(t, \xi)d\xi -x_k(t) \Big|\le \frac {4\gt}{|g(0)|} \wt \vep.   
\end{equation*}
  }    

The meaning of \lq\lq{}$\vep$ is small enough\rq\rq{} in Theorem B is that $\wt \vep$ defined by $\wt \vep={\vep}^{1/3}$ satisfies some inequalities like \eqref{cond00}. 
Most theorems on the WSST and FSST analysis are stated in the form of Theorem B, see e.g.  \cite{Daub_Lu_Wu11,Thakur_Wu11,Wu_thesis, MOM14,BMO18}. Actualy the statements of part (b) and part (c)  in Theorem A give us more direct bounds of the estimates. 
We call the quantity on the left-hand side (LHS) of \eqref{transformation_est0} the IF estimate error, and  call that on LHS of \eqref{reconstr_1st0} the error of component recovery (or component separation). 
The statements in Theorem A can be found in \cite{Wu_thesis, MOM14, BMO18} but with some different IF estimate errors. 
For example, \cite{MOM14,BMO18} gave IF estimate error 
$\frac {\vep}{\wt \vep}\big(\Gamma_0(t) (\gt+2\phi_k^\gp(t)) +\frac 1{2\pi} \wt \Gamma_0(t)\big)$ instead of 
$\frac {\vep}{\wt \vep}\big(\Gamma_0(t) \gt +\frac 1{2\pi} \wt \Gamma_0(t)\big)$ in \eqref{transformation_est0}. 
One can also find that Theorem A is a special case of Theorem 1 in Section 3 below (see Remark \ref{rem:back_regular_FSST}).


\bigskip 

Observe that the condition \eqref{cond_phi_2nd_der} or \eqref{cond_phi_2nd_der1}
requires the slow change of the IF $\phi^\gp_k(t)$ of each component $x_k(t)$. There is no mathematical guarantee for the IF estimate and the component separation for a multicomponent signal $x(t)$ with a component $x_k(t)$ having a fast-changing frequency (e.g. $\phi_k^{\gp\gp}(t)$, the changing rate of  IF of $x_k(t)$, is not very small).  To this regard, the 2nd-order FSST was introduced in \cite{MOM15} and later the 2nd-order WSST was proposed in \cite{OM17} with the theoretical analysis of  the 2nd-order FSST established in \cite{BMO18}. 

Suppose $V_x(t, \eta)\not=0$ and $\partial_t\big(\frac{\partial _\eta V_x(t, \eta)}{V_x(t, \eta)}\big)\not=i 2\pi $. Denote 
$$
\wt q(t, \eta)=\frac {\partial_t\Big(\frac{\partial _t V_x(t, \eta)}{V_x(t, \eta)}\Big)}{i2\pi-\partial_t\Big(\frac{\partial _\eta V_x(t, \eta)}{V_x(t, \eta)}\Big)}. 
$$
The 2nd-order FSST in \cite{BMO18} is defined as 
\begin{equation*}
R_{x, \gga}^{2nd, \gl}(t, \xi)
=\int_{\{|V_x(t, \eta)|>\gga}  V_x(t, \eta) \frac 1{\gl}h\big(\frac{\xi-\go^{2nd}_x(t, \eta)}\gl\big) d\eta,  
\end{equation*}
where $\go_x^{2nd}(t, \eta)$ is the phase transformation for the 2nd-order FSST: for $(t, \eta)$ with $V_x(t, \eta)\not=0$,  
\begin{equation*}
\go_x^{2nd}(t, \eta)=\left\{
\begin{array}{ll}
{\rm Re}\big\{\frac{\frac {\partial}{\partial t} V_x(t, \eta)}{2\pi i V_x(t, \eta)}\big\}
+{\rm Re}\big\{ \wt q(t, \eta)  \partial_t\big(\frac{\partial _\eta V_x(t, \eta)}{V_x(t, \eta)}\big)  \big\},  &\hbox{if $\partial_t\big(\frac{\partial _\eta V_x(t, \eta)}{V_x(t, \eta)}\big)\not=i 2\pi $,}\\
{\rm Re}\big\{\frac{\frac {\partial}{\partial t} V_x(t, \eta)}{2\pi i V_x(t, \eta)}\big\}, &\hbox{elsewhere.}
\end{array}
\right. 
\end{equation*}

Let $\vep>0$ and $\gt>0$. ${B}^{(2)}_{\vep, \gt}$ denote the set of multicomponent signals of \eqref{AHM} 
satisfying 
\eqref{cond_basic}, the well-separated condition \eqref{freq_resolution}, and the following conditions: 
\begin{eqnarray}
 \label{cond_basic_2nd} && A_k(t)\in C^2(\R)\cap L_\infty(\R), \phi_k(t)\in C^3(\R),  
\phi^{\gp\gp}_k(t) \in L_\infty(\R), \\
\label{cond_phi_3nd_der}&& |A^\gp_k(t)|\le \vep, \; |A_k^{\gp\gp}(t)|\le \vep,  \;  
|\phi^{(3)}_k(t)|\le \vep, \; t\in \R. 
\end{eqnarray}
Then for $x(t)\in {B}^{(2)}_{\vep, \gt}$,  statements for the 2nd-order FSST similar to those in Theorem B hold, under some conditions which are more complicated than \eqref{cond00} 
 since  no band-limited 
 window functions and the 2nd-order  phase transformation $\go_x^{2nd}(t, \eta)$ are involved. See  \cite{BMO18} for the details. Observe that  there is no direct  boundedness restriction on 
$\phi_k^{\gp\gp}(t)$ in the definition of ${B}^{(2)}_{\vep, \gt}$. 


\subsection{Adaptive FSST with a time-varying parameter}

We consider the window function given by
\begin{equation}
\label{window_general}
g_\gs(t)=\frac 1\gs g(\frac t\gs),
\end{equation}
where $\gs>0$ is a parameter, $g(t)$ is a function in $L_2(\R)$ with $g(0)\not=0$ and having certain decaying order as $t\rightarrow \infty$. If
\begin{equation}
\label{def_g}
g(t)=\frac 1{\sqrt {2\pi}} e^{-\frac {t^2}2},
\end{equation}
then $g_\gs(t)$ is the Gaussian window function. The parameter $\gs$ is also called the window width in the time-domain of the window function $g_\gs(t)$ since the time duration  $\Delta_{g_\gs}$ of $g_\gs$ is $\gs$ (up to a constant): $\Delta_{g_\gs}=\gs\Delta_{g}$, 
where $\Delta_{g}$ is the time duration of $g$.

For a signal $x(t)$, the STFT of $x(t)$ with a time-varying parameter is defined in \cite{LCHJJ18} as
\begin{eqnarray}
\label{def_STFT_para1}\label{def_STFT_para2}
\wt V_{x}(t, \eta) 
\hskip -0.6cm &&=\int_\R x(\tau)g_{\gs(t)}(\tau-t)e^{-i2\pi \eta (\tau-t)}d\tau
=\int_\R x(t+\tau)\frac 1{\gs(t)}g (\frac \tau{\gs(t)})e^{-i2\pi \eta\tau}d\tau, 
\end{eqnarray}
where $\gs=\gs(t)$ is a positive function of $t$.  
$\wt V_{x}(t, \eta)$ is called the adaptive STFT of $x(t)$ with $g_\gs$. 

\bigskip 

Before we move on to review the SST associated with the adaptive STFT, we introduce some notations used in this and next sections. Denote 
$$
g_1(\tau)=\tau g(\tau), \; g_2(\tau )=\tau ^2g(\tau ), \; g_3(\tau)=\tau g^\gp(\tau ). 
$$
Thus 
$$
g_{1, \gs}(\tau)=\frac {\tau}{\gs^2} g(\frac \tau\gs), \; 
g_{2, \gs}(\tau)=\frac {\tau^2}{\gs^3} g(\frac \tau\gs), \;
g_{3, \gs}(\tau)=\frac {\tau}{\gs^2} g^\gp(\frac \tau\gs). 
$$
We use $\wt V^{g_j}_x(t, \eta)$ and $\wt V^{g^\gp}_x(t, \eta)$ to denote the adaptive STFT defined by \eqref{def_STFT_para1} with $g_\gs$ replaced by $g_{j, \gs}$ and $g^\gp_\gs(\tau)=\frac 1{\gs}g^\gp(\frac \tau\gs)$ respectively, where $1\le j\le 3$. 

\bigskip

For  $x(t)=A  e^{i2\pi c t}$, one can show that (see \cite{LCHJJ18})  
if $\wt V_x(t, \eta)\not=0$, then 
\begin{equation}
\label{def_phase_para_complex}
\go^{adp, c}_x(t, \eta)=\frac {\frac{\partial}{\partial t} \wt V_x(t, \eta)}{i2\pi \wt V_x(t, \eta)}+\frac {\gs'(t)}{i2\pi \gs(t)}
+ \frac {\gs'(t)}{\gs(t)}\frac {\wt V^{g_3}_x(t, \eta)}{i2\pi \wt V_x(t, \eta)}, 
\end{equation}
is $c$, the IF of $x(t)$. 
Hence, for a general $x(t)$, at $(t, \eta)$,  \cite{LCHJJ18} defines 
the real part of the quantity of $\go^{adp, c}_x(t, \eta)$ in the above equation, denoted by $\go^{adp}_x(t, \eta)$, 
as the phase transformation of the adaptive FSST: 
\begin{equation*}
\go^{adp}_x(t, \eta)={\rm Re}\Big\{\frac{\partial_t\big(\wt V_x(t, \eta)\big)}{i2\pi \wt V_x(t, \eta)}\Big\}+ \frac {\gs'(t)}{\gs(t)} {\rm Re}\Big\{ \frac {\wt V^{g_3}_x(t, \eta)}{i2\pi \wt V_x(t, \eta)} \Big\},  \quad \hbox{for $\wt V_x(t, \eta)\not=0$}.
\end{equation*}
Then the (1st-order) adaptive FSST, denoted by 
 $R_{x, \gga}^{adp, \gl}$, is defined by 
\begin{equation}
\label{def_adpFSST}
R_{x, \gga}^{adp, \gl}(t, \xi)
=\int_{|\wt V_x(t, \eta)|>\gga} \wt V_x(t, \eta) \frac 1{\gl}h\big(\frac{\xi-\go^{adp}_x(t, \eta)}\gl\big) d\eta, 
\end{equation}
where $\gga>0, \gl>0$ and $h(t)$ is a compactly supported function as describled in \S2.1. 

\bigskip


Next we consider the 2nd-order adaptive FSST.  For a linear chirp signal, 
\begin{equation}
\label{def_chip_At}
x(t)=A e^{i2\pi \phi(t)}=A e^{i2\pi (ct +\frac 12 r t^2)}, 
\end{equation}
it was shown in \cite{LCHJJ18} that $\go_{x}^{adp, 2nd, c}$ defined below is 
$c+rt$, the IF of $x(t)$: 
\begin{equation}
\label{def_transformation_2nd_complex}
\go_{x}^{adp, 2nd, c}=\frac{\gs'(t)}{i2\pi \gs(t)}+ \frac {\frac{\partial}{\partial t} \wt V_x(t, \eta)}{i2\pi \wt V_x(t, \eta)}
- \frac{\wt V^{g_1}_x(t, \eta)}{i2\pi \wt V_x(t, \eta)} P_0(t, \eta)
+ \frac {\gs'(t)}{\gs(t)}\frac {\wt V^{g_3}_x(t, \eta)}{i2\pi \wt V_x(t, \eta)}, 
\end{equation}
for $(t, \eta)$ satisfying $\frac{\partial}{\partial \eta}\Big( \frac {\wt V^{g_1}_x(t, \eta)}{\wt V_x(t, \eta)}\Big)\not=0$ and $\wt V_x(t, \eta)\not=0$, where 
 \begin{equation}
\label{def_P0}
P_0(t, \eta)=\frac 1{\frac{\partial}{\partial \eta}\Big( \frac {\wt V^{g_1}_x(t, \eta)}{\wt V_x(t, \eta)}\Big) }\Big\{\frac{\partial}{\partial \eta}\Big(\frac {\frac{\partial}{\partial t} \wt V_x(t, \eta)}{\wt V_x(t, \eta)}\Big)+ \frac {\gs'(t)}{\gs(t)}
\frac{\partial}{\partial \eta}\Big(\frac {\wt V^{g_3}_x(t, \eta)}{\wt V_x(t, \eta)}\Big)\Big\}.
\end{equation}
Thus  the authors of \cite{LCHJJ18} define the real part of $\go_{x}^{adp, 2nd, c}$ as the phase transformation for the 2nd-order adaptive FSST. Namely,  the phase transformation $\go_{x}^{adp, 2nd}$ is defined by 
\begin{equation*}
\go^{adp, 2nd}_x(t, \eta)=\left\{
\begin{array}{l}
{\rm Re}\Big\{\frac {\frac{\partial}{\partial t} \wt V_x(t, \eta)}{i2\pi \wt V_x(t, \eta)}\Big\}
- {\rm Re}\Big\{ \frac{\wt V^{g_1}_x(t, \eta)}{i2\pi \wt V_x(t, \eta)} P_0(t, \eta)\Big\}+ \frac {\gs'(t)}{\gs(t)} {\rm Re}\Big\{ \frac {\wt V^{g_3}_x(t, \eta)}{i2\pi \wt V_x(t, \eta)} \Big\},\\
\hskip 5cm \hbox{if $\frac{\partial}{\partial \eta}\Big(\frac {\wt V^{g_1}_x(t, \eta)}{\wt V_x(t, \eta)}\Big)\not=0$ and $\wt V_x(t, \eta)\not=0;$}\\
{\rm Re}\Big\{\frac {\frac{\partial}{\partial t} \wt V_x(t, \eta)}{i2\pi \wt V_x(t, \eta)}\Big\}+ \frac {\gs'(t)}{\gs(t)} {\rm Re}\Big\{ \frac {\wt V^{g_3}_x(t, \eta)}{i2\pi \wt V_x(t, \eta)} \Big\},
\hbox{if $\frac{\partial}{\partial \eta}\Big(\frac {\wt V^{g_1}_x(t, \eta)}{\wt V_x(t, \eta)}\Big)=0$,  $\wt V_x(t, \eta)\not=0$,}
\end{array}
\right.
\end{equation*}

Finally we define the 2nd-order adaptive FSST with $V_x(t, \eta)\not =0$ and 
$\frac{\partial}{\partial \eta}\Big(\frac {\wt V^{g_1}_x(t, \eta)}{\wt V_x(t, \eta)}\Big)\not=0$ described by  thresholds $\gga_1>0, \gga_2>0$. More precisely, 
define 
$$
\go_{x, \gga_1, \gga_2}^{apd, 2nd, c}(t, \eta)=\left\{
\begin{array}{l}
\hbox{\rm quantity in \eqref{def_transformation_2nd_complex},  \quad  if  $|V_x(t, \eta)|>\gga_1$ and $\big|\frac{\partial}{\partial \eta}\Big(\frac {\wt V^{g_1}_x(t, \eta)}{\wt V_x(t, \eta)}\big| $} >\gga_2, \\
\hbox{\rm quantity in \eqref{def_phase_para_complex},  \quad  if  $|V_x(t, \eta)|>\gga_1$ and $\big|\frac{\partial}{\partial \eta}\Big(\frac {\wt V^{g_1}_x(t, \eta)}{\wt V_x(t, \eta)}\big| $} \le \gga_2. 
\end{array}
\right. 
$$
Let $\go_{x, \gga_1, \gga_2}^{apd, 2nd}(t, \eta)$=Re$\big(\go_{x, \gga_1, \gga_2}^{apd, 2nd, c}(t, \eta)\big)$. 

Again, let $h(t)$ be a compactly supported function with certain smoothness and $\int_\R h(t)dt=1$. We define the 2nd-order adaptive FSST $R_{x, \gga_1, \gga_2}^{adp, 2nd,\gl}$ by 
\begin{equation}
\label{def_adp2ndFSST}
R_{x, \gga_1, \gga_2}^{adp, 2nd,\gl}(t, \xi)
=\int_{\big\{\eta: \; |\wt V_x(t, \eta)|>\gga_1, \; \big|{\partial \eta}\big({\wt V^{g_1}_x(t, \eta)}/{\wt V_x(t, \eta)}\big)\big| >\gga_2\big\}} \wt V_x(t, \eta) \frac 1{\gl}h\big(\frac{\xi-\go^{adp, 2nd}_{x,\gga_1, \gga_2}(t, \eta)}\gl\big) d\eta. 
\end{equation}


\section{Analysis of adaptive FSST}

We assume 
\begin{equation}
\label{freq_resolution_adp}
d\rq{}=\min_{k\in \{2, \cdots, K\}}\min_{t\in \R}(\phi^\gp_k(t)-\phi^\gp_{k-1}(t))> 0.
\end{equation}
Thus $x(t)$ satisfies the well-separated condition  \eqref{freq_resolution} with resolution $=d\rq{}/2$. 
However, the value $d\rq{}$ may be very small. In this case, we cannot apply Theorem A directly. 
The reason is that to guarantee the results in Theorem A to hold, the window function $g$ needs to satisfy supp$(\wh g)\subseteq [-\frac {d\rq{}}2, \frac {d\rq{}}2]$. 
If $d\rq{}$ is quite small, then $g$ has a very good frequency resolution, which implies by the uncertainty principle that  $g$  has a very poor time resolution, or equivalently $g$ has a very large time duration, 
which results in large errors in the IF estimate and component recovery 
(see Remark \ref{rem:window_width}).  
We use  the adaptive STFT and FSST to adjust the time-varying  window width $\gs(t)$ at certain local time $t$ where 
 the frequencies of two components are close.

In this section we consider the case that each component $x_k(t)=A_k(t)e^{i2\pi \phi_k(t)}$ is approximated locally by a sinusoidal signal. More precisely, we assume 
$A_k^\gp(t)$ and $\phi_k^{\gp\gp}(t)$  are small:
\begin{equation}
\label{condition1} |A_k^\gp(t)|\le \vep_1, \; |\phi_k^{\gp\gp}(t)|\le \vep_2, \; t\in \R, \; 1\le k\le K, 
\end{equation}
for some positive number $\vep_1, \vep_2$. Let ${\cal D}_{\vep_1, \vep_2}$ denote the set of multicomponent signals of \eqref{AHM} satisfying \eqref{cond_basic0}, \eqref{cond_basic},  \eqref{freq_resolution_adp} and \eqref{condition1}. 

\bigskip 

Let $x(t)\in {\cal D}_{\vep_1, \vep_2}$.   Write $x_k(t+\tau)$ as 
\begin{eqnarray*}
x_k(t+\tau)\hskip -0.6cm &&=x_k(t)e^{i2\pi \phi^\gp_k(t)\tau }+(A_k(t+\tau)-A_k(t))e^{i2\pi \phi_k(t+\tau)}\\
&&\qquad +x_k(t)e^{i2\pi \phi^\gp_k(t)\tau}\big(
e^{i2\pi (\phi_k(t+\tau)-\phi_k(t)-\phi_k^\gp(t) \tau)}-1\big).
\end{eqnarray*}
Then we have 
\begin{eqnarray}
\nonumber
\wt V_{x}(t, \eta)\hskip -0.6cm &&=\sum_{k=1}^K \int_\R x_k(t+\tau)\frac 1{\gs(t)}g(\frac \tau{\gs(t)}) e^{-i2\pi \eta \tau}d\tau\\
\nonumber 
&&=\sum_{k=1}^K \int_\R x_k(t)e^{i2\pi \phi_k^\gp(t) \tau}\frac 1{\gs(t)}g(\frac \tau{\gs(t)}) e^{-i2\pi \eta \tau}d\tau +\rem_0, 
\end{eqnarray}
or 
\begin{equation}
\label{STFT_approx_1st}
\wt V_{x}(t, \eta)
=\sum_{k=1}^K x_k(t)  \wh g\big(\gs(t)(\eta-\phi_k^\gp(t) \big) +\rem_0,  
\end{equation}
where $\rem_0$ is the remainder for the expansion of  $\wt V_{x}(t, \eta)$  
in \eqref{STFT_approx_1st} given by 
\begin{eqnarray}
\label{def_rem0}
&&\rem_0=\sum_{k=1}^K \int_\R \Big\{ 
(A_k(t+\tau)-A_k(t))e^{i2\pi \phi_k(t+\tau)}
\\
\nonumber  &&\qquad \qquad 
+x_k(t)e^{i2\pi \phi^\gp_k(t)\tau}\big(
e^{i2\pi (\phi_k(t+\tau)-\phi_k(t)-\phi_k^\gp(t) \tau)}-1\big)
\Big\}\frac 1{\gs(t)}g(\frac \tau{\gs(t)}) e^{-i2\pi \eta \tau}d\tau. 
\end{eqnarray}
With $|A_k(t+\tau)-A_k(t)|\le  \vep_1|\tau|$ 
and 
$$
|e^{i2\pi (\phi_k(t+\tau)-\phi_k(t)-\phi_k^\gp(t) \tau)}-1|\le 2\pi |\phi_k(t+\tau)-\phi_k(t)-\phi_k^\gp(t) \tau|
\le 
\pi\vep_2 |\tau|^2, 
$$
we have  
\begin{eqnarray*}
|\rem_0|\hskip -0.6cm && \le \sum_{k=1}^K \int_\R  \vep_1 | \tau| \frac 1{\gs(t)}|g(\frac \tau{\gs(t)})| d\tau+ \sum_{k=1}^K A_k(t)\int_\R  \pi \vep_2 | \tau|^2 \frac 1{\gs(t)}|g(\frac \tau{\gs(t)})|
 d\tau   \\
 &&=K \vep_1 I_1 \gs(t)  +\pi  \vep_2 I_2 \gs^2(t) \sum_{k=1}^K A_k(t), 
\end{eqnarray*}
where $I_n$ is defined in \eqref{def_In}. 
Hence we have 
\begin{equation}
\label{rem0_est}
|\rem_0|\le \gs(t) \Lambda_0(t),  
\end{equation}
where 
\begin{equation}
\label{def_Lam0}
\Lambda_0(t)=K \vep_1  I_1 +\pi  \vep_2 I_2 \gs (t) \sum_{k=1}^K A_k(t). 
\end{equation}

 $\wt V^{g^\gp}_x(t, \eta)$ 
can be expanded as \eqref{STFT_approx_1st}
with remainder
 $\rem_0^\gp$, 
defined as $\rem_0$ in \eqref{def_rem0} with $g(\tau)$ replaced by 
 $g^\gp(\tau)$. 
Then we have the estimate for the remainder similar to \eqref{rem0_est}. More precisely, we have 
\begin{equation}
\label{err_est}
|\rem_0^\gp|\le \gs(t) \wt \Lambda_0(t), 
\end{equation}
where 
\begin{equation}
\label{def_tLam0}
\wt \Lambda_0(t)=K \vep_1  \wt I_1 +\pi \vep_2 \wt I_2 \gs(t) \sum_{k=1}^K A_k(t), 
\end{equation}
with $\wt I_n$ defined in \eqref{def_tIn}.

\bigskip 
\begin{mrem}
\label{rem:ConditionB}
Condition \eqref{condition1} is essentially the condition  \eqref{cond_phi_2nd_der1}. If $A_k(t), \phi_k(t)$ satisfy  \eqref{cond_phi_2nd_der}, then we have a similar error bound for the expansion of $\wt V_x(t, \eta)$. More precisely, suppose  $A_k(t), \phi_k(t)$ satisfy 
\begin{equation}
\label{conditionB}
  |A^\gp_k(t)|\le \vep_1 \phi^\gp_k(t),  \;  |\phi^{\gp\gp}_k(t)|\le \vep_2\phi^\gp_k(t), \; t\in \R, \; 
 M^{\gp\gp}_k=\sup_{t\in \R} |\phi^{\gp\gp}_k(t)|<\infty. 
\end{equation}
Then {\rm (}see {\rm \cite{Daub_Lu_Wu11})} 
\begin{eqnarray*}
&&|A_k(t+\tau)-A_k(t)|\le  \vep_1|\tau| ( \phi_k^\gp(t)+\frac 12 M_k^{\gp\gp} |\tau|), \\
&&|\phi_k(t+\tau)-\phi_k(t)-\phi_k^\gp(t) \tau| \le \vep_2 \tau^2 ( \frac 12 \phi_k^\gp(t)+\frac 16 M_k^{\gp\gp} |\tau|). 
\end{eqnarray*}
Thus, 
we can expand $\wt V_x(t, \eta)$ as \eqref{STFT_approx_1st} with $|\rem_0|\le \gs(t) \Lambda_0(t)$, where in this case  $\Lambda_0(t)$ is 
\begin{equation}
\label{CondtionB_est}
\Lambda_0(t)=
\vep_1 \sum_{k=1}^K \big(I_1 \phi^\gp_k(t)+\frac 12 M_k^{\gp\gp} I_2\gs(t)\big)
   +\pi  \vep_2 \gs(t) \sum_{k=1}^K A_k(t)\big(I_2\phi^\gp_k(t)+\frac 13 M_k^{\gp\gp} I_3\gs(t)\big). 
\end{equation}
With the condition of \eqref{conditionB}, we have an estimate $\gs(t) \wt \Lambda_0(t)$ for $\rem_0\rq{}$ with 
$\wt \Lambda_0(t)$ defined by \eqref{CondtionB_est} with $I_j$ replaced by $\wt I_j$. 
In this paper we consider the condition \eqref{condition1}. 
The statements for theoretical analysis of the adaptive FSST with condition \eqref{conditionB} instead of  \eqref{condition1}
 are still valid as long as $\Lambda_0(t)$ 
in \eqref{def_Lam0}, $\Lambda\rq{}_0(t)$ in \eqref{def_tLam0} and so on are replaced respectively by that in \eqref{CondtionB_est} and similar terms. This also applies to the 2nd-order adaptive FSST in Section 4, where we will not repeat again this discussion on the condition like \eqref{conditionB}. 
\hfill $\blacksquare$ 
\end{mrem}

If the remainder $\rem_0$ in \eqref{STFT_approx_1st} is small, then 
the term $x_k(t)\wh g\big(\gs(t)(\eta-\phi_k^\gp(t) \big)$  in \eqref{STFT_approx_1st} givens the time-frequency zone of the STFT $\wt V_{x_k}$ of the $k$th component $x_k(t)$ of $x(t)$. If  in addition, $g$ is band-limited, that is 
 $\wh g$ is compactly supported, to say supp($\wh g)\subset [-\gt, \gt]$, then $x_k(t)\wh g\big(\gs(t)(\eta-\phi_k^\gp(t) \big)$ lies  within the zone:
 $$
\{(t, \eta): |\eta-\phi_k^\gp(t)|< \frac {\gt}{\gs(t)}, t\in \R\}. 
 $$
Thus the multicomponent signal $x(t)$ is well-separated (that is $Z_k\cap Z_{\ell}=\O, k\not=\ell$.), provided that $\gs(t)$ satisfies 
\begin{equation}
\label{separated_cond_1st_compact}
\gs(t)\ge \frac {2\gt}{ \phi_k^\gp(t)-\phi_{k-1}^\gp(t)}, \; t\in \R, k=2, \cdots, K. 
\end{equation}
Observe that our well-separated condition \eqref{separated_cond_1st_compact} is different from that in \eqref{freq_resolution} considered in \cite{Wu_thesis} and \cite{MOM14}. 

If $\wh g$ is not compactly supported,  
we need to define the ``support'' of  $\wh g$ outside which $\wh g(\xi)\approx 0$. 
More precisely, for a given threshold $0<\tau_0<1$, if
$|\wh g(\xi)|\le \tau_0$ for $|\xi|\ge \ga$, then we say $\wh g(\xi)$ is ``supported\rq\rq{} in $[-\ga, \ga]$. 
When $|\wh g(\xi)|$ is even and decreasing for $\xi\ge 0$,   
then $\ga$ can be obtained by solving 
\begin{equation}
\label{def_ga_general}
|\wh g(\ga)|=\tau_0. 
\end{equation}
For example, when  $g$ is the Gaussian function given by \eqref{def_g}, 
then, with $\wh g(\xi)=e^{-2\pi^2 \xi^2}$,
the corresponding $\ga$ is given by  
\begin{equation}
\label{def_ga}
\ga=\frac 1{2\pi}\sqrt{2\ln (\frac 1{\tau_0 })} . 
\end{equation}

For $g$ with $\wh g(\xi)$ ``supported\rq\rq{} in $[-\ga, \ga]$,  we then define the time-frequency zone $Z_k$ of the $k$th-component $x_k(t)$ of $x(t)$ by 
\begin{equation}
\label{def_Zk}
Z_k=\{(t, \eta): |\wh g\big(\gs(t)(\eta-\phi_k^\gp(t) \big) |>\tau_0, t\in \R\}=
\{(t, \eta): |\eta-\phi_k^\gp(t)|< \frac {\ga}{\gs(t)}, t\in \R\}. 
\end{equation}
Thus the multicomponent signal $x(t)$ is well-separated, if $\gs(t)$ satisfies 
\begin{equation}
\label{separated_cond_1st}
\gs(t)\ge \frac {2\ga}{ \phi_k^\gp(t)-\phi_{k-1}^\gp(t)}, \; t\in \R, k=2, \cdots, K. 
\end{equation}
In this case $Z_k\cap Z_{\ell}=\O, k\not=\ell$. 
In this section we assume that \eqref{separated_cond_1st} holds for some $\gs(t)$. Due to \eqref{freq_resolution_adp}, there always exists $C^1(\R)$ and bounded  $\gs(t)$ such that \eqref{separated_cond_1st} holds. 
For the sinusoidal function-based adaptive FSST, the authors in \cite{LCHJJ18} suggest to choose $\gs(t)$ 
as 
 $$
 \gs_1(t)=\max\{\frac {2\ga}{ \phi_k^\gp(t)-\phi_{k-1}^\gp(t)}, k=2, \cdots, K\}.
 $$

 Observe that for  $\gs(t)$ satisfying \eqref{separated_cond_1st}, since $\phi_k^\gp(t)$ is bounded, we know 
$$
\|\frac 1{\gs(t)}\|_\infty<\infty. 
$$
In addition,  in this case 
\begin{equation}
\label{different_IF_est}
\gs(t)| \phi_k^\gp(t)-\phi_\ell^\gp(t)|\ge{2\ga |k-\ell|}.  
\end{equation}
 
\bigskip

Next we will present our analysis results on the adaptive FSST in Theorem \ref{theo:main_1st} below, 
where $\ga$ is defined by \eqref{def_ga_general}, and 
$\sum_{\ell\not=k}$ denotes $\sum_{\ell\in \{1, \cdots, K\}\backslash\{k\}}$. 
$\wt V_x(t, \eta)$ is the adaptive STFT of $x(t)$
with such a  window function  $g\in \cS$ that $|\wh g(\xi)|$ is  even and decreasing for $\xi\ge 0$. 


\begin{theo}
\label{theo:main_1st} Let $x(t)\in {\cal D}_{\vep_1, \vep_2}$ for some small $\vep_1, \vep_2>0$. Then we have the following.

{\rm (a)} Suppose $\wt \ep_1$ satisfies  $\wt \ep_1 \ge \gs(t) \Lambda_0(t)+ \tau_0 \sum_{k=1}^K A_k(t)$. 
Then for $(t, \eta)$ with $|\wt V_x(t, \eta)|>\wt \ep_1$, there exists $k\in \{1, 2, \cdots, K\}$ such that $(t, \eta)\in Z_k$.   

{\rm (b)} For $(t, \eta)$ with $|\wt V_x(t, \eta)|\not=0$, we have  
\begin{equation}
\label{transformation_approx}
\go_x^{apd, c}(t, \eta)-\phi^\gp_k(t)= \frac{\Rem_1}{i2\pi \wt V_x(t, \eta)},  
\end{equation} 
where  
\begin{eqnarray}
\label{def_Rem1}
&&\Rem_1=i2\pi \big(\eta-\phi_k^\gp(t)\big)\rem_0 -  \frac{ \rem_0^\gp}{\gs(t)}+
i2\pi \sum_{\ell\not=k} x_\ell(t) \big(\phi^\gp_\ell(t)-\phi^\gp_k(t)\big)
\wh g\big(\gs(t)(\eta-\phi_\ell^\gp(t)))\big).
 \end{eqnarray}
Hence,  for $(t, \eta)$ satisfying $|\wt V_x(t, \eta)|>\wt \ep_1$ and  $(t, \eta)\in Z_k$, we have 
\begin{equation}
\label{transformation_est}
|\go_x^{apd}(t, \eta)-\phi^\gp_k(t)|<  \bd_1, 
\end{equation}
where 
\begin{eqnarray}
\nonumber 
\bd_1\hskip -0.6cm &&=
\frac 1{\wt \ep_1}\big(\ga \Lambda_0(t)+\frac 1{2\pi}\wt \Lambda_0(t)\big)\\
\label{def_bd1}&&  \quad +\frac 1{\wt \ep_1} \max_{k\in\{1, \cdots, K\}}\Big\{ \sum_{\ell\not=k} A_\ell(t)|\phi_\ell^\gp (t)-\phi_k^\gp (t)| \sup_{\{u: |u|<\ga\} }
\big|\wh g\big(u+\gs(t)(\phi_k^\gp(t)-\phi_\ell^\gp(t))\big)\big|\Big\}. 
\end{eqnarray}

{\rm (c)} Suppose that $\wt \ep_1$ satisfies the condition in part {\rm (a)} 
and that $bd_1$ in part {\rm (b)} satisfies $bd_1\le \frac \ga{\gs(t)}$. Then for $\wt\ep_3$ satisfying
 $bd_1\le \wt \ep_3\le \frac \ga{\gs(t)}$, we have  
\begin{equation}
 \label{reconstr_1st}
 \Big| \lim_{\gl\to 0} \frac{\gs(t)}{g(0)}\int_{|\xi-\phi^\gp_k(t)|<\wt \ep_3} R_{x, \wt \ep_1}^{adp,  \gl}(t, \xi)d\xi -x_k(t) \Big|\le \bd_2, 
\end{equation}
where 
\begin{equation}
 \label{def_bd2}
\bd_2=\frac 1{|g(0)|} \Big\{ 2\ga (\gs(t)\Lambda_0(t)+\wt \ep_1)+A_k(t)\; \big|\int_{|u|\ge \ga}  \wh g(u) du \big|+\sum_{\ell\not=k }A_\ell(t) m_{\ell, k}(t)\Big\}, 
\end{equation}
with
$$
m_{\ell, k}(t)
=  \big|\int_{|u|\le \ga }   \wh g\big(u+\gs(t)(\phi_k^\gp(t)-\phi_\ell^\gp(t))\big) du \big|. 
$$ 
\end{theo}

\bigskip 
\begin{mrem}
\label{rem:bandlimited}
When $\wh g(\xi)$ is supported in $[-\ga, \ga]$, we can set $\tau_0$ in Theorem \ref{theo:main_1st} part {\rm (a)} to be zero, and thus,  the condition for part  {\rm (a)} is for $\wt \ep_1$ to satisfy  $\wt \ep_1 \ge \gs(t) \Lambda_0(t)$. In addition, in this case the 2nd term in \eqref{def_bd1} is zero, and  
$m_{\ell, k}(t)$ and $\int_{|u|\ge \ga}  \wh g(u) du$ in \eqref{def_bd2} are also zero. 
 Hence in this case,  \eqref{def_bd1} and \eqref{reconstr_1st} are respectively reduced to 
 \begin{equation}
 \label{IF_est_AFSST_compact}
 |\go_x^{apd}(t, \eta)-\phi^\gp_k(t)|\le \frac 1{\wt \ep_1}\big(\ga \Lambda_0(t)+\frac 1{2\pi}\wt \Lambda_0(t)\big), 
 \end{equation}
 and 
 \begin{equation}
 \label{component_est_AFSST_compact}
 \Big| \lim_{\gl\to 0} \frac{\gs(t)}{g(0)}\int_{|\xi-\phi^\gp_k(t)|<\wt \ep_3} R_{x, \wt \ep_1}^{adp,  \gl}(t, \xi)d\xi -x_k(t) \Big|\le \frac {2\ga (\gs(t)\Lambda_0(t)+\wt \ep_1)}{|g(0)|}, 
 \end{equation}
 for any  $\wt \ep_3$ with 
 $\frac 1{\wt \ep_1}\big(\ga \Lambda_0(t)+\frac 1{2\pi}\wt \Lambda_0(t)\big) \le \wt \ep_3\le \frac {\ga}{\gs(t)}$. 
 
Furthermore, the statement of Theorem \ref{theo:main_1st} can be written in the form of Theorem B. 
 For simplicity, we just consider the case $\vep_1=\vep_2$. Write    
 $\Lambda_0(t), \wt \Lambda_0(t)$ defined by  \eqref{def_Lam0} and \eqref{def_tLam0} respectively as 
 $$
 \Lambda_0(t)=\vep_1\gl_0(t), \; \wt \Lambda_0(t)=\vep_1\wt \gl_0(t),
$$
with 
$$
 \gl_0(t)= K  I_1 +\pi   I_2 \gs (t) \sum_{k=1}^K A_k(t), \; 
 \wt \gl_0(t)= K \wt I_1 +\pi  \wt I_2 \gs (t) \sum_{k=1}^K A_k(t).  
 $$
 Let $\wt \ep_1=\vep_1^3$. If $\vep_1$ is small enough such that 
 \begin{equation}
 \label{ep1_cond}
 \wt \ep_1\le \max\{\ga \; \|\frac1{\gs(t)}\|_\infty, \|\ga \lambda_0(t)+\frac 1{2\pi}\wt \lambda_0(t)\|^{-1}_\infty\}, 
 \end{equation}
  then   
  $$
  \frac 1{\wt \ep_1}\big(\ga \Lambda_0(t)+\frac 1{2\pi}\wt \Lambda_0(t)\big) \le {\wt \ep_1}, \; \hbox{and 
  $\wt \ep_1 \le \frac{\ga}{\gs(t)}$}. 
  $$
  Thus we have the following corollary.  
  \hfill $\blacksquare$
  
 \end{mrem}

\begin{cor}
\label{cor:theo_standand_form}  Suppose $x(t)\in {\cal D}_{\vep_1, \vep_1}$ for some small $\vep_1>0$, and  supp($\wh g) \subseteq [-\ga, \ga]$. Let $\wt \ep_1=\vep_1^3$. If $\vep_1$ is small enough such that \eqref{ep1_cond} holds, then we have the following. 
  
  {\rm (a)} 
 For $(t, \eta)$ satisfying  $|\wt V_x(t, \eta)|>\wt \ep_1$,  there exists a unique $k\in \{1, 2, \cdots, K\}$ such that $(t, \eta)\in Z_k$.   

{\rm (b)} Suppose $(t, \eta)$ satisfies $|\wt V_x(t, \eta)|>\wt \ep_1$ and  $(t, \eta)\in Z_k$. Then 
\begin{equation*}
|\go_x^{apd}(t, \eta)-\phi^\gp_k(t)|<  \wt \ep_1.  
\end{equation*}

{\rm (c)} 
For any $k$, $1\le k\le K$, 
\begin{equation*}
 \Big| \lim_{\gl\to 0} \frac{\gs(t)}{g(0)}\int_{|\xi-\phi^\gp_k(t)|<\wt \ep_1} R_{x, \wt \ep_1}^{adp,  \gl}(t, \xi)d\xi -x_k(t) \Big|\le \frac{2\ga (\gs(t)\Lambda_0(t)+\wt \ep_1)}{|g(0)|}.  
\end{equation*}

\end{cor}
\begin{mrem}
 \label{rem:back_regular_FSST}
When $\gs(t)\equiv 1$, $R_{x, \wt \ep_1}^{adp,  \gl}(t, \xi)$ is the regular FSST 
$R_{x, \wt \ep_1}^{\gl}(t, \xi)$ defined by \eqref{def_FSST}. Suppose supp$(\wh g)\subseteq [-\ga, \ga]$. 
Then \eqref{transformation_est} and \eqref{reconstr_1st} {\rm (}which, by Remark \ref{rem:bandlimited}, are 
\eqref{IF_est_AFSST_compact} and \eqref{component_est_AFSST_compact}{\rm )} with $\vep_1=\vep_2$  are 
 respectively \eqref{transformation_est0} and \eqref{reconstr_1st0} with $\vep=\vep_1, \gt=\ga$. Thus in the case $\gs(t)\equiv 1$, 
 Theorem \ref{theo:main_1st} is reduced to Theorem A, and Corollary \ref{cor:theo_standand_form} is Theorem B. 
\hfill $\blacksquare$
\end{mrem}

 \begin{mrem}
 \label{rem:window_width}
 Observe that $\Lambda_0(t)$ and $\wt \Lambda_0(t)$ defined by \eqref{def_Lam0} and \eqref{def_tLam0} respectively depend on $\gs(t)$. 
Smaller $\gs(t)$ results in smaller $\Lambda_0(t)$ and $\wt \Lambda_0(t)$, hence  
the corresponding $\go_x^{apd}(t, \eta)$ provides a more accurate  estimate for $\phi^\gp_k(t)$ as implied by \eqref{IF_est_AFSST_compact} and the adaptive FSST in \eqref{reconstr_1st} gives a recovery of $x_k(t)$ with a smaller error as shown in \eqref{component_est_AFSST_compact}. 
\end{mrem}

\begin{mrem}
\label{rem:theo_standand_form}
When $\wh g(\xi)$ is not supported on $[-\ga, \ga]$, but $|\wh g(\xi)|$  decays fast as $|\xi|\to \infty$, then the terms in the summation for $\bd_1$ in \eqref{def_bd1} will be small as long as $\tau_0$ is small.  Recall that we assume $|\wh g(\xi)|$ is even and decreasing.
Then for $\ell=k-1$, since $\gs(t)(\phi_k^\gp(t)-\phi_{k-1}^\gp(t))> 2\ga$,  we have for $|u|\le \ga$, 
 $$
 \big| \wh g \big(u+\gs(t)(\phi_k^\gp(t)-\phi_{k-1}^\gp(t))\big)\big|\le | \wh g  \big(u+2\ga\big)| \le   |\wh g(\ga)| = \tau_0
$$
Similarly, we have $\sup_{\{|u|\le \ga\}} \big| \wh g \big(u+\gs(t)(\phi_k^\gp(t)-\phi_{k+1}^\gp(t))\big)\big| \le \tau_0$. 
The quantities $\sup_{\{|u|\le \ga\}} \big| \wh g \big(u+\gs(t)(\phi_k^\gp(t)-\phi_\ell^\gp(t))\big)\big|$ for other $\ell \not= k-1, k, k+1$ are much smaller than $\tau_0$ since 
$\gs(t)(\phi_k^\gp(t)-\phi_\ell^\gp(t))$ is large {\rm (}see \eqref{different_IF_est}{\rm )} and $\wh g$ is rapidly decreasing. Thus the summation in $\bd_1$ 
is dominated by  $\tau_0 
\max_{k=2, \cdots K-1} \{A_{k+1}(t)(\phi_{k+1}^\gp (t)-\phi_k^\gp (t)),  A_{k-1}(t)(\phi_{k}^\gp (t)-\phi_{k-1}^\gp (t))\}$.  

\bigskip 

The functions $m_{\ell, k}(t)$ in \eqref{def_bd2} could be small if $\tau_0$ is small. More precisely, for $\ell=k-1$, 
we have 
\begin{eqnarray*}
&&m_{k-1, k}(t)\le  \int_{|u|\le  \ga } | \wh g \big( u+\gs(t)(\phi_k^\gp(t)-\phi_{k-1}^\gp(t) ) \big)| du\\
&& \le \int_{|u|\le \ga } | \wh g  \big(u+2\ga\big)| du \le  2 \ga  |\wh g(\ga)| =2\ga \tau_0. 
\end{eqnarray*}
We can show similarly that $m_{k+1, k}(t)\le 2\ga \tau_0$. For other $\ell$, $m_{\ell, k}(t)$ will be much smaller than $2\ga \tau_0$.

Finally, let us loot at the term $\big|\int_{|u|\ge \ga}  \wh g(u) du \big|$. It will be small if $\ga$ is large. Here we give an estimate of this term when $g(t)$ is the Gaussian function defined by \eqref{def_g}. In this case, one can obtain
 
\begin{equation}
\label{est_Guasian_intgral}
 \int_{|u|\ge \ga
 } e^{-2\pi^2 u^2}du\le \frac 1{\sqrt {2\pi}} \frac{e^{-2\pi^2 \ga^2}}{1+\sqrt{1-e^{-2\pi^2 \ga^2}}}=  \frac 1{\sqrt {2\pi}} \frac{\tau_0}{1+\sqrt{1-\tau_0}}. 
\end{equation}
To summarize, for $g$ with $\wh g(\xi)$ decaying rapidly as $\xi\to \infty$,  the statements in Corollary \ref{cor:theo_standand_form} hold if the same conditions are assumed and that  $\tau_0$ is smaller enough (and hence $\ga$ is large enough). 
\hfill $\blacksquare$
\end{mrem}

\bigskip 

In the rest of section, we give the proof of Theorem \ref{theo:main_1st}. 

\bigskip 

{\bf Proof  of  Theorem \ref{theo:main_1st} Part (a)}.  Assume $(t, \eta)\not \in \cup _{k=1}^K Z_k$. Then for any $k$,  by the definition of $Z_k$ in \eqref{def_Zk}, we have  
$|\wh g\big(\gs(t)(\eta-\phi_k^\gp(t))\big)|\le \tau_0$. Hence, by \eqref{STFT_approx_1st} and \eqref{rem0_est}, we have 
\begin{eqnarray*}
|\wt V_x(t, \eta)|\hskip -0.6cm &&\le \sum_{k=1}^K |x_k(t) \wh g\big(\gs(t)(\eta-\phi_k^\gp(t))\big)| + |\rem_0|
\\
&&\le \tau _0 \sum_{k=1}^K A_k(t)  +\gs(t) \Lambda_0(t)\le \wt \ep_1, 
\end{eqnarray*}
a contradiction to the assumption $|\wt V_x(t, \eta)|>\wt \ep_1$. Hence the statement in (a) holds. 
\hfill $\blacksquare$

\bigskip 

{\bf Proof  of  Theorem \ref{theo:main_1st} Part (b)}.   By a direct calculation, we have 
\begin{equation}
\label{result_partial_tV}
\pd_t \wt V_x(t, \eta)=\big(i2\pi \eta-\frac {\gs'(t)}{\gs(t)}\big)\wt V_x(t, \eta)
- \frac {\gs'(t)}{\gs(t)}\wt V^{g_3}_x(t, \eta)-\frac 1{\gs(t)} \wt V^{g^\gp}_x(t, \eta). 
\end{equation}
By \eqref{STFT_approx_1st} with $g$ replaced by $g^\gp$, 
\begin{eqnarray*}
&&\wt V_x^{g^\gp}(t,\eta)=\sum_{\ell=1}^K \int_\R x_\ell(t)e^{i2\pi  \phi_\ell^\gp(t) \tau}\frac 1{\gs(t)}g^\gp(\frac \tau{\gs(t)}) e^{-i2\pi \eta \tau}d\tau +\rem_0^\gp\\
&&=\sum_{\ell=1}^K x_\ell(t)(g^\gp)^{\wedge}\big(\gs(t)(\eta-\phi_\ell^\gp(t))\big)
 +\rem_0^\gp\\
&&=i2\pi \gs(t)\sum_{\ell=1}^K x_\ell(t)(\eta-\phi_\ell^\gp (t))\wh g\big(\gs(t)(\eta-\phi_\ell^\gp(t))\big)
 +\rem_0^\gp. 
\end{eqnarray*}
This and \eqref{result_partial_tV} imply that 
\begin{eqnarray*}
&& \big( \go_x^{apd, c}(t, \eta)-\phi^\gp_k(t)\big)i2\pi \wt V_x(t,\eta)\\
&&=\pd_t \wt V_x(t,\eta)+\frac {\gs'(t)}{\gs(t)}\big(\wt V_x(t,\eta)+\wt V_x^{g_3}(t,\eta)\big)
-i2\pi \phi^\gp_k(t)\wt V_x(t,\eta)\\
&&=i2\pi \eta \wt V_x(t,\eta)-\frac 1{\gs(t)}\wt V_x^{g^\gp}(t,\eta)-i2\pi \phi^\gp_k(t)\wt V_x(t,\eta)\\
&&=i2\pi \big(\eta -\phi_k^\gp(t)\big)\Big (\sum_{\ell=1}^K x_\ell(t)
\wh g\big(\gs(t)(\eta-\phi_\ell^\gp(t))\big)+\rem_0\Big)\\
&&\qquad -\frac 1{\gs(t)}\Big(
i2\pi \gs(t)\sum_{\ell=1}^K x_\ell(t)(\eta-\phi_\ell^\gp (t))\wh g\big(\gs(t)(\eta-\phi_\ell^\gp(t))\big)
 +\rem_0^\gp
\Big)\\
&&=i2\pi \big(\eta -\phi_k^\gp(t)\big)\rem_0-\frac{\rem_0^\gp}{\gs(t)}
+i2\pi \sum_{\ell\not=k} x_\ell(t)(\phi_\ell^\gp (t)-\phi_k^\gp (t))\wh g\big(\gs(t)(\eta-\phi_\ell^\gp(t))\big)\\
&&=\Rem_1.  
\end{eqnarray*}
This shows \eqref{transformation_approx}.

For \eqref{transformation_est}, with the assumptions   
$|\wt V_x(t, \eta)|>\wt \ep_1$,  
we have 
\begin{eqnarray*}
&&|\go_x^{apd}(t, \eta)-\phi^\gp_k(t)|\le |\go_x^{apd, c}(t, \eta)-\phi^\gp_k(t)|\\ 
&&\le \Big| \frac{\Rem_1}{i2\pi \wt V_x(t, \eta)}
\Big|< \frac{|\Rem_1|}{2\pi \wt \ep_1}\\  
&&\le \frac 1{\wt \ep_1}\Big\{ |\big(\eta -\phi_k^\gp(t)\big)| \gs(t)\Lambda_0(t)+\frac 1{2\pi}\wt \Lambda_0(t)+\sum_{\ell\not=k} A_\ell(t)|\phi_\ell^\gp (t)-\phi_k^\gp (t)| \; 
|\wh g\big(\gs(t)(\eta-\phi_\ell^\gp(t))\big)|\Big\}\\
&&\le \frac 1{\wt \ep_1}\Big\{\ga \Lambda_0(t)+\frac 1{2\pi}\wt \Lambda_0(t)+ \sum_{\ell\not=k} A_\ell(t)|\phi_\ell^\gp (t)-\phi_k^\gp (t)| \; 
|\wh g\big(\gs(t)(\eta-\phi_\ell^\gp(t))\big)|\Big\}\\
&&\le \bd_1,   
\end{eqnarray*}
as desired, where the second last inequality follows from $|\eta -\phi_k^\gp(t)|<\frac\ga{\gs(t)}$ since $(t, \eta)\in Z_k$. 
\hfill $\blacksquare$

\bigskip 

{\bf Proof  of  Theorem \ref{theo:main_1st} Part (c)}. First we have the following result which can be derived as that on p.254 in \cite{Daub_Lu_Wu11}:
\begin{equation}
\label{FSST_STFT_relation1_1st}
 \lim_{\gl\to 0}\int_{|\xi-\phi^\gp_k(t)|<\wt \ep_3} R_{x, \wt \ep_1}^{adp, \gl}(t, \xi)d\xi =
 \int_{X_t}  \wt V_x(t, \eta) 
 d\eta ,
\end{equation}
where 
$$
X_t=\big\{\eta:\; |\wt V_x(t, \eta)|>\wt \ep_1 \; 
\hbox{and $\big|\phi^\gp_k(t)-\go_x^{adp}(t, \eta)\big|<\wt \ep_3$}\big \}. 
$$

Let 
$$
Y_t=\big\{\eta:\; |\wt V_x(t, \eta)|>\wt \ep_1 \; \hbox{and $(t, \eta)\in Z_k$}\big \}. 
$$
By Theorem \ref{theo:main_1st} part (b), if $\eta\in Y_t$, then $\big|\phi^\gp_k(t)-\go_x^{adp}(t, \eta)\big|< \bd_1\le \wt \ep_3$. Thus $\eta\in X_t$. Hence $Y_t\subseteq X_t$. 

On the other hand, suppose $\eta\in X_t$. Since $|\wt V_x(t, \eta)|>\wt \ep_1$, by Theorem \ref{theo:main_1st} part (a), $(t, \eta)\in Z_\ell$ for an $\ell$ in $\{1, 2, \cdots, K\}$. If $\ell\not =k$, then by Theorem Theorem \ref{theo:main_1st} part (b),
\begin{eqnarray*}
\big|\phi^\gp_k(t)-\go_x^{adp}(t, \eta)\big|\hskip -0.6cm &&\ge |\phi^\gp_k(t) -\phi^\gp_\ell(t)|-\big|\phi^\gp_\ell(t)-\go_x^{adp}(t, \eta)\big|\\
&& > \frac{2\ga}{\gs(t)}-\bd_1\ge \frac{2\ga}{\gs(t)}-\wt \ep_3\ge \wt \ep_3, 
\end{eqnarray*}
since $\bd_1\le \wt \ep_3\le \frac{\ga}{\gs(t)}$. 
This contradictes to the assumption $\big|\phi^\gp_k(t)-\go_x^{adp}(t, \eta)\big|<\wt \ep_3$ since $\eta \in X_t$. Hence $\ell=k$ and $\eta\in Y_t$. Thus we get $X_t=Y_t$. Therefore, from \eqref{FSST_STFT_relation1_1st}, we have 
\begin{equation}
\label{FSST_STFT_relation2_1st}
 \lim_{\gl\to 0}\int_{|\xi-\phi^\gp_k(t)|<\wt \ep_3} R_{x, \wt \ep_1}^{adp,\gl}(t, \xi)d\xi =
 \int_{\{|\wt V_x(t, \eta)|>\wt \ep_1\}\cap \{\eta: (t, \eta)\in Z_k\}}  \wt V_x(t, \eta) 
 d\eta. 
\end{equation}

Furthermore, 
\begin{eqnarray*}
&&\big| \int_{\{|\wt V_x(t, \eta)|>\wt \ep_1\}\cap \{\eta: (t, \eta)\in Z_k\}}  \wt V_x(t, \eta) d\eta 
-\frac{g(0)}{\gs(t)}x_k(t)\big|\\
&&= \big| \int_{\{\eta: (t, \eta)\in Z_k\}}  \wt V_x(t, \eta) d\eta 
-\frac{g(0)}{\gs(t)}x_k(t)- \int_{\{|\wt V_x(t, \eta)|\le\wt \ep_1\}\cap \{\eta: (t, \eta)\in Z_k\}}  \wt V_x(t, \eta) d\eta 
\big|\\
&&\le  \wt \ep_1 \frac{2\ga}{\gs(t)} +\big| \int_{\{\eta: (t, \eta)\in Z_k\}}  \big(\sum_{\ell=1}^K x_\ell(t) 
\wh g\big(\gs(t)(\eta-\phi_\ell^\gp(t))\big)
+\rem_0 \big)d\eta  
-\frac{g(0)}{\gs(t)}x_k(t)\big| \\
&& \le \wt \ep_1 \frac{2\ga}{\gs(t)} +|\rem_0| \; \frac{2\ga}{\gs(t)}+ 
\big| \int_{|\eta-\phi_k^\gp(t)|<\frac \ga{\gs(t)}}   x_k(t) \wh g\big(\gs(t)(\eta-\phi_k^\gp(t))\big) d\eta  
-\frac{g(0)}{\gs(t)}x_k(t)\big| \\
&&\qquad + \sum_{\ell\not=k }A_\ell(t) |\int_{|\eta-\phi_k^\gp(t)|<\frac \ga{\gs(t)}}\wh g\big(\gs(t)(\eta-\phi_\ell^\gp(t))\big)  d\eta| \\
&&= (|\rem_0|+\wt \ep_1) \; \frac{2\ga}{\gs(t)} +
\big|\frac {x_k(t)}{\gs(t)} \int_{|u|< \ga }   \wh g(u) du  
-\frac{g(0)}{\gs(t)}x_k(t)\big| \\
&&\qquad +\sum_{\ell\not=k }\frac{A_\ell(t)}{\gs(t)} 
 \big|\int_{|u|<\ga }   \wh g\big(u+\gs(t)(\phi_k^\gp(t)-\phi_\ell^\gp(t))\big) du \big|\\
&&=(|\rem_0|+\wt \ep_1)  \; \frac{2\ga}{\gs(t)}
\big|\frac {x_k(t)}{\gs(t)}  \int_\R   \wh g(u) du  
-\frac{g(0)}{\gs(t)}x_k(t)-\frac{x_k(t)}{\gs(t)}\int_{|u|\ge \ga}
 \wh g(u) du \big|
 +\sum_{\ell\not=k } \frac{A_\ell(t)}{\gs(t)}  m_{\ell, k}(t) \\
&&= (|\rem_0|+\wt \ep_1) \; \frac{2\ga}{\gs(t)}
+\big|x_k(t) \frac{g(0)}{\gs(t)}
-\frac{g(0)}{\gs(t)}x_k(t)-\frac{x_k(t)}{\gs(t)} \int_{|u|\ge \ga}
 \wh g(u) du \big|+\frac 1{\gs(t)}\sum_{\ell\not=k } A_\ell(t)  m_{\ell, k}(t) \\
&&\le   
\big(\gs(t)\Lambda_0(t)+\wt \ep_1\big) \; \frac{2\ga}{\gs(t)}
+\frac{A_k(t)}{\gs(t)}\big|\int_{|u|\ge \ga}
 \wh g(u) du \big|+\frac 1{\gs(t)} \sum_{\ell\not=k }A_\ell(t) m_{\ell, k}(t). 
\end{eqnarray*}
The above estimate, together with \eqref{FSST_STFT_relation2_1st},
leads to \eqref{reconstr_1st}. 
This completes the proof of Theorem \ref{theo:main_1st} Part (c). 
\hfill $\blacksquare$

\section{Analysis of 2nd-order adaptive FSST}


We consider multicomponent signals $x(t)$ of \eqref{AHM} 
satisfying \eqref{cond_basic_2nd} and being well approximated locally by a linear chirp signal of \eqref{def_chip_At}
with 
$A_k^\gp(t)$ and $\phi_k^{(3)}(t)$ are small:
\begin{equation}
\label{condition2}
|A_k^\gp(t)|\le \vep_1, \; |\phi_k^{(3)}(t)|\le \vep_3, \; t\in \R, \; 1\le k\le K, 
\end{equation}
for some positive number $\vep_1, \vep_3$,

For a given $t$, we use $G_k(\xi)$ to denote the Fourier transform of $e^{i\pi \gs(t)\phi^{\gp\gp}_k(t) \tau^2}g(\tau)$, namely, 
$$
G_k(\xi)={\cal F}\Big(e^{i\pi \gs(t)\phi^{\gp\gp}_k(t) \tau^2}g(\tau)\Big)\big(\xi)=
\int_{\R} e^{i\pi\gs(t) \phi^{\gp\gp}_k(t) \tau^2}g(\tau) e^{-i2\pi \xi \tau}d\tau,   
$$
where ${\cal F}$ denotes the Fourier transform. 
Note that $G_k(\xi)$ depends on $t$ also if $\phi^{\gp\gp}_k(t)\not=0$. We drop $t$ in $G_\ell$ 
for simplicity.

For each component $x_k(t)=A_k(t)e^{i2\pi \phi_k(t)}$, we write $x_k(t+\tau)$ as 
\begin{eqnarray*}
x_k(t+\tau)\hskip -0.6cm &&=x_k(t)e^{i2\pi (\phi^\gp_k(t)\tau+\frac 12\phi^{\gp \gp}_k(t)\tau^2) }
+(A_k(t+\tau)-A_k(t))e^{i2\pi \phi_k(t+\tau)}\\
&&\quad +x_k(t)e^{i2\pi (\phi^\gp_k(t)\tau+\frac 12\phi^{\gp \gp}_k(t)\tau^2) }
\big(
e^{i2\pi (\phi_k(t+\tau)-\phi_k(t)-\phi_k^\gp(t) \tau- \frac 12\phi^{\gp \gp}_k(t)\tau^2)}-1\big).
\end{eqnarray*}
Then we have 
\begin{eqnarray}
\nonumber
\wt V_{x}(t, \eta)\hskip -0.6cm &&=\sum_{k=1}^K \int_\R x_k(t+\tau)\frac 1{\gs(t)}g(\frac \tau{\gs(t)}) e^{-i2\pi \eta \tau}d\tau\\
\label{STFT_approx0}
&&=\sum_{k=1}^K \int_\R x_k(t)e^{i2\pi (\phi_k^\gp(t) \tau +\frac12\phi^{\gp\gp}_k(t) \tau^2)}\frac 1{\gs(t)}g(\frac \tau{\gs(t)}) e^{-i2\pi \eta \tau}d\tau +\err_0\\
\label{STFT_approx}
&&=\sum_{k=1}^K x_k(t)  
G_k\big(\gs(t)(\eta-\phi_k^\gp(t) )\big)+\err_0, 
\end{eqnarray}
where 
\begin{eqnarray}
\label{def_err0}
&&\err_0=\sum_{k=1}^K \int_\R \Big\{ 
(A_k(t+\tau)-A_k(t))e^{i2\pi \phi_k(t+\tau)}
\\
\nonumber  &&
+x_k(t)e^{i2\pi (\phi^\gp_k(t)\tau+\frac 12\phi^{\gp \gp}_k(t)\tau^2) }
\big(
e^{i2\pi (\phi_k(t+\tau)-\phi_k(t)-\phi_k^\gp(t) \tau- \frac 12\phi^{\gp \gp}_k(t)\tau^2)}-1\big)
\Big\}\frac 1{\gs(t)}g(\frac \tau{\gs(t)}) e^{-i2\pi \eta \tau}d\tau. 
\end{eqnarray}

To distinguish the different types of the remainders for the expansion of  $\wt V_{x}(t, \eta)$ resulted from different local approximations for $x_k(t+\tau)$,   
in this section we use \lq\lq{}$\err$\rq\rq{}, which means residual,  to denote  the remainder  for the expansion of  $\wt V_{x}(t, \eta)$  in \eqref{STFT_approx0}. 

With $|A_k(t+\tau)-A_k(t)|\le \vep_1|\tau|$ 
and 
$$
|e^{i2\pi (\phi_k(t+\tau)-\phi_k(t)-\phi_k^\gp(t) \tau- \frac 12\phi^{\gp \gp}_k(t)\tau^2)}-1|\le 
2\pi \frac 16 \sup_{\eta\in \R}|\phi^{(3)}_k(\eta) \tau^3|
\le \frac \pi 3 \vep_2 |\tau|^3, 
$$
we have  
\begin{eqnarray*}
&&|\err_0|\le   \sum_{k=1}^K  \int_\R  \vep_1 | \tau| \frac 1{\gs(t)}|g(\frac \tau{\gs(t)})| d\tau+
\sum_{k=1}^K A_k(t)\int_\R  \frac \pi3 \vep_3 | \tau|^3 \frac 1{\gs(t)}|g(\frac \tau{\gs(t)})|
 d\tau \\
 &&=K \vep_1 I_1 \gs(t)  +\frac \pi 3 \vep_3 I_3 \gs^3(t) \sum_{k=1}^K A_k(t), 
\end{eqnarray*}
where $I_n$ is defined in \eqref{def_In}. 
Hence we have 
\begin{equation}
\label{err0_est}
|\err_0|\le \gs(t) \Pi_0(t),  
\end{equation}
where 
$$
\Pi_0(t)=K \vep_1  I_1 +\frac \pi 3 \vep_3 I_3 \gs^2(t) \sum_{k=1}^K A_k(t). 
$$

Thus if $\vep_1, \vep_2$ are small enough, then $|\err_0|$ is small and hence, $G_k\big(\gs(t)(\eta-\phi_k^\gp(t) )\big)$ provides the time-frequency zone for $\wt V_{x_k}(t, \eta)$.  In the following we  describe those  time-frequency zones mathematically. 
Let $0<\tau_0<1$ be a given small number as the threshold. 
Denote 
\begin{equation}
\label{def_Ok}
O_k=
\{(t, \eta): |G_k\big(\gs(t)(\eta-\phi_k^\gp)\big)|>\tau_0, t\in \R\}. 
\end{equation}
We assume  again $|G_k(\xi)|$ is even and decreasing for $\xi\ge 0$. 
Then $O_k$ can be written as 
\begin{equation}
\label{def_Ok2}
O_k=\{(t, \eta): |\eta-\phi_k^\gp(t)|< \frac {\ga_k}{\gs(t)}, t\in \R\}. 
\end{equation}
where $\ga_k=\ga_k(t)$ is obtain by solving 
$|G_k(\xi)|=\tau_0. 
$
We will assume the multicomponent signal $x(t)$ is well-separated, that is there is $\gs(t)$ such that 
\begin{equation}
\label{cond_no_overlapping} 
O_k\cap O_{\ell}=\O, \quad k\not=\ell. 
\end{equation}
  
As an example, let us look at what are $O_k$ and $\ga_k(t)$ look like when $g$ is the Gaussian function defined by  \eqref{def_g}. 
One can obtain for this $g$ (see \cite{LCHJJ18}),  
\begin{equation}
\label{STFT_LinearChip}
G_{k}(u)=\frac {1}{\sqrt{1-i2\pi \phi_k^{\gp\gp}(t)\gs^2(t)}}\;
e^{-\frac{2\pi^2 u^2}{1+(2\pi \phi_k^{\gp\gp}(t)\gs^2(t))^2} (1+i2\pi \phi_k^{\gp\gp}(t)\gs^2(t))}. 
\end{equation}
Thus 
\begin{equation}
\label{abs_Gk}
|G_k(u)|=\frac 1{\big(1+(2\pi \phi_k^{\gp\gp}(t)\gs^2(t))^2\big)^{\frac 14}}\;
e^{-\frac{2\pi^2}{1+(2\pi \phi_k^{\gp\gp}(t)\gs^2(t))^2}u^2}. 
\end{equation} 
Therefore, in this case, assuming $\tau_0 (1+  (2\pi \phi^{\gp\gp}_k(t) \gs^2(t))^2)^{\frac 14}\le 1$, 
\begin{equation}
\label{def_gak_Gaussian}
\ga_k=\sqrt{1+ 
 (2\pi \phi^{\gp\gp}_k(t) \gs^2(t))^2} \; \frac 1{2\pi}\sqrt{2\ln (\frac 1{\tau_0})-\frac 12 \ln(1+ 
 (2\pi \phi^{\gp\gp}_k(t) \gs^2(t))^2)}. 
\end{equation}
Authors of \cite{LCHJJ18} consider  a larger zone $O^\gp_k$ in the time-frequency plane:
   \begin{equation}
\label{def_zone_larger}
 O^\gp_k=\Big\{(t, \eta): |\eta-\phi_k^\gp(t)|< \frac {\ga}{\gs(t)}
  \Big(1+ 
 2\pi |\phi^{\gp\gp}_k(t)| \gs^2(t)\Big), t\in \R\Big\}.   
\end{equation}
where $\ga$ is defined by \eqref{def_ga}. 
They obtain that if for 
$k=2, \cdots, K, $
\begin{eqnarray}
\label{LinearChirp_sep_cond1}
&& 4\ga \sqrt{\pi} \sqrt {|\phi''_k(t)|+|\phi''_{k-1}(t)|}\le \phi'_k(t)-\phi'_{k-1}(t), \hbox{and} \\
\label{LinearChirp_sep_cond2}
 &&  \max_{2\le k\le K}\Big\{
 \frac{4\ga}{b_k(t)+\sqrt{b_k(t)^2-8\ga a_k(t)}}
 \Big\}\le \min_{2\le k\le K}\Big\{
\frac{4\ga}{b_k(t)-\sqrt{b_k(t)^2-8\ga a_k(t)}}
\Big\},
 \end{eqnarray}
then the components $x_k(t), 1\le k\le K$ of $x(t)$ are well-separated in the time-frequency plane in the sense that $O_k\rq{}\cap O_{\ell}\rq{}=\O,  k\not=\ell$, where 
\begin{equation*}
\label{def_akbk}
a_k(t)=2\pi\ga (|\phi''_{k-1}(t)| +|\phi''_k(t)|),  \; b_k(t)=\phi_{k}'(t)-\phi_{k-1}'(t).
\end{equation*}
 \cite{LCHJJ18} shows that any $\gs(t)$ between the two quantities in \eqref{LinearChirp_sep_cond2} separates the components $x_k(t)$ of $x(t)$ in the time-frequency plane,  and 
it suggests to choose $\gs(t)$ to be 
\begin{equation}
\label{def_gs2}
\gs_2(t)=
\max_{2\le k\le K}\Big\{\frac{4\ga}{b_k(t)+\sqrt{b_k(t)^2-8\ga a_k(t)}} \Big\}. 
 \end{equation}




Let $g$ be a window function in the Schwarz class $\cS$ with  $|\wh g(\xi)|$ is  even and decreasing for $\xi\ge 0$. Let ${\cal D}^{(2)}_{\vep_1, \vep_2}$ denote the set of multicomponent signals of \eqref{AHM} satisfying 
\eqref{cond_basic_2nd}, \eqref{freq_resolution_adp}, \eqref{condition2} and  that $x(t)$ is well-sperated with $g$, that is there is $\gs(t)$ such that \eqref{cond_no_overlapping} holds. 


\bigskip 

We introduce more notations to describe our main theorem on the 2nd-order adptive FSST. For $j\ge 0$, denote 
\begin{eqnarray}
\label{def_Gjk}
G_{j,k}(t, \eta)\hskip -0.6cm &&=\int_{\R} e^{i2\pi (\phi_k^\gp(t) \tau +\frac12\phi^{\gp\gp}_k(t) \tau^2)}\frac {\tau^j}{\gs(t)^{j+1}}g(\frac \tau{\gs(t)}) e^{-i2\pi \eta \tau}d\tau\\
\nonumber &&={\cal F}\Big(e^{i\pi \gs(t)\phi^{\gp\gp}_k(t) \tau^2}\tau^j g(\tau)\Big)\big(\gs(t)(\eta-\phi_k^\gp(t))\big).   
\end{eqnarray}
Clearly, 
$$
G_{0, k}(t, \eta)=G_k\big(\gs(t)(\eta-\phi_k^\gp(t) )\big), 
$$
and one can obtain for $j\ge 1$, 
\begin{eqnarray}
&&\label{GjkwithGk}
G_{j, k}(t, \eta)=\frac 1{(-i2\pi)^j}G^{(j)}_k\big(\gs(t)(\eta-\phi_k^\gp(t))\big). 
\end{eqnarray}
 
Let $\err_1$, $\err_2$, $\err_0^\gp$, and $\err_1^\gp$ be the residuals defined as $\err_0$ in \eqref{def_err0} with $g(\tau)$ replaced respectively by $g_1(\tau)$,  $g_2(\tau)$,  $g^\gp(\tau)$, and $g_3(\tau)=\tau g^\gp(\tau)$. 
Then we have the estimates for these residuals similar to \eqref{err0_est}. More precisely, we have 
\begin{equation}
\label{err_est}
|\err_1|\le \gs(t) \Pi_1(t), \;  |\err_2|\le \gs(t) \Pi_2(t), \quad |\err_0^\gp|\le \gs(t) \wt \Pi_0(t), \; 
|\err_1^\gp|\le \gs(t) \wt \Pi_1(t),      
\end{equation}
where 
\begin{eqnarray*}
&&\Pi_1(t)=K \vep_1  I_2 +\frac \pi 3 \vep_3 I_4 \gs^2(t) \sum_{k=1}^K A_k(t), \;  \Pi_2(t)=K \vep_1  I_3 +\frac \pi 3 \vep_3 I_5 \gs^2(t) \sum_{k=1}^K A_k(t), \\ 
&&\wt \Pi_0(t)=K \vep_1  \wt I_1 +\frac \pi 3 \vep_3 \wt I_3 \gs^2(t) \sum_{k=1}^K A_k(t), \;  
\wt \Pi_1(t)=K \vep_1  \wt I_2 +\frac \pi 3 \vep_3 \wt I_4 \gs^2(t) \sum_{k=1}^K A_k(t), 
\end{eqnarray*}
with $I_n, \wt I_n$ defined in \eqref{def_tIn}. 

Denote 
\begin{eqnarray*}
\label{def_Bk}
&&B_k(t, \eta)=\sum_{\ell\not=k} x_\ell(t) \big(\phi^\gp_\ell(t)-\phi^\gp_k(t)\big)G_{0, \ell}(t, \eta), 
\quad D_k(t, \eta)=\sum_{\ell\not=k} x_\ell(t) \big(\phi^{\gp\gp}_\ell(t)-\phi^{\gp\gp}_k(t)\big)G_{1, \ell}(t, \eta),\\
\label{def_Ek}
&&E_k(t, \eta)=\sum_{\ell\not=k} x_\ell(t) \big(\phi^\gp_\ell(t)-\phi^\gp_k(t)\big)G_{1, \ell}(t, \eta),
\quad F_k(t, \eta)=\sum_{\ell\not=k} x_\ell(t) \big(\phi^{\gp\gp}_\ell(t)-\phi^{\gp\gp}_k(t)\big)G_{2, \ell}(t, \eta), 
\end{eqnarray*}
and 
\begin{eqnarray}
\label{def_Err1}
\Err_1\hskip -0.6 cm &&=i2\pi B_k(t, \eta)+i 2\pi \gs(t)D_k(t, \eta)
+i2\pi \big(\eta-\phi_k^\gp(t)\big)\err_0 -  \frac{ \err_0^\gp}{\gs(t)}-i2\pi \phi_k^{\gp\gp}(t) \gs(t)\err_1,  \\
\label{def_Err2}
\Err_2\hskip -0.6 cm &&=4\pi^2\gs(t) E_k(t, \eta)+4\pi^2 \gs^2(t)F_k(t, \eta)
\\
\nonumber && \qquad 
+i 2\pi  \; \err_0+4\pi^2 \big(\eta-\phi_k^\gp(t)\big)\gs(t)\; \err_1 +  i2\pi \; 
\err_1^\gp-4\pi^2 \phi_k^{\gp\gp}(t) \gs^2(t)\; \err_2.
 \end{eqnarray}

Next we will provide Theorem \ref{theo:main} on the 2nd-order adaptive FSST, which also consists of parts (a)-(c). The proof of part (b) is based on the following two lemmas whose proof is postponed to Appendix. 
\begin{lem}
\label{lem:lem1} Let $\Err_1$ be the quantity definde by \eqref{def_Err1}. Then 
\begin{equation}
\label{result_lem1}
\pd_t \wt V_x(t, \eta)=\big(i2\pi \phi^\gp_k(t)-\frac {\gs'(t)}{\gs(t)}\big)\wt V_x(t, \eta)+
i2\pi \phi_k^{\prime\prime}(t)\gs(t)\wt V^{g_1}_x(t, \eta)- \frac {\gs'(t)}{\gs(t)}\wt V^{g_3}_x(t, \eta)+\Err_1. 
\end{equation}
\end{lem}


\begin{lem}
\label{lem:lem3} For $(t, \eta)$ satisfying 
$V_x(t, \eta)\not =0$ and 
$\frac{\partial}{\partial \eta}\Big(\frac {\wt V^{g_1}_x(t, \eta)}{\wt V_x(t, \eta)}\Big)\not=0$, we have  
\begin{equation}
\label{result_lem3}
P_0(t, \eta)=i2\pi \gs(t)\phi^{\prime\prime}_k(t)+\Err_3,  
\end{equation}
where 
\begin{equation}
\label{Err3}
\Err_3=\frac {\wt V_x(t, \eta) \; \Err_2-\pd_\eta \wt V_x(t, \eta) \; \Err_1 }{\wt V_x(t, \eta)\pd_\eta 
\wt V^{g_1}_x(t, \eta)-\wt V^{g_1}_x(t, \eta)\pd_\eta 
\wt V_x(t, \eta)}, 
\end{equation}
with $\Err_1$ and $\Err_2$ defined by \eqref{def_Err1} and \eqref{def_Err2} respectively. 
\end{lem}

\bigskip 
Denote 
\begin{equation}
\label{def_Mlk}
M_{\ell, k}(t)=\gs(t)\int_{\{\eta: \; (t, \eta)\in O_k\}}| G_{0, \ell}(t, \eta)| d\eta
= \int_{|u|< \ga_k } | G_\ell \big(u+\gs(t)(\phi_k^\gp(t)-\phi_\ell^\gp(t))\big) | du.   
\end{equation}

\begin{theo}
\label{theo:main} Suppose $x(t)\in {\cal D}^{(2)}_{\vep_1, \vep_2}$ with a window function $g(t)\in {\cal S}$ for some small $\vep_1, \vep_2>0$. Then we have the following. 

{\rm (a)} Suppose $\wt \vep_1$ satisfies  $\wt \vep_1 \ge \vep_0 \sum_{k=1}^K A_k(t)+ \gs(t) \Pi_0(t)$. Then for $(t, \eta)$ with $|\wt V_x(t, \eta)|>\wt \vep_1$, there exists $k\in \{1, 2, \cdots, K\}$ such that $(t, \eta)\in O_k$.   

${\rm (b)}$ Suppose $(t, \eta)$ satisfies $|\wt V_x(t, \eta)|>\wt \vep_1$,  $|\pd _\eta \big(\wt V^{g_1}_x(t, \eta)/\wt V_x(t, \eta)\big)|>\wt \vep_2$,  and $(t, \eta)\in O_k$. 
Then 
\begin{equation}
\label{transformation_approx_2nd}
\go_x^{apd, 2nd, c}(t, \eta)-\phi^\gp_k(t)= \Err_4, 
\end{equation} 
where 
\begin{equation*}
\Err_4=\frac{\Err_1}{i2\pi \wt V_x(t, \eta)}- \frac{\wt V^{g_1}_x(t, \eta)\Err_3}{i2\pi\wt V_x(t, \eta)}.   
\end{equation*}
Furthermore,  
\begin{equation}
\label{IF_est_Bd1}
|\go_x^{apd, 2nd}(t, \eta)-\phi^\gp_k(t)|<  \Bd_1, 
\end{equation}
where 
\begin{equation}
\label{def_Bd1}
\Bd_1=\max_{1\le k\le K}\sup_{\eta\in O_k}\Big\{\frac{|\Err_1|}{2\pi \wt \vep_1}+ \frac1{2\pi \wt \vep_1^3 \wt \vep_2}
|\wt V^{g_1}_x(t, \eta)| \big(|\pd_\eta \wt V_x(t, \eta)|\;  |\Err_1| + \wt \vep_1 |\Err_2|\big)\Big\}.    
\end{equation}

{\rm (c)} Suppose that $\wt \vep_1$ satisfies the condition in part {\rm (a)} and  $\Bd_1 \le \frac 12 L_k(t)$, 
where 
\begin{equation}
\label{defLk}
L_k(t)=\frac 1{\gs(t)}\min\{\ga_k+\ga_{k-1}, \ga_k+\ga_{k+1}\}.  
\end{equation}
Then for any $\wt \vep_3=\wt \vep_3(t)>0$ satisfying $\Bd_1\le \wt \vep_3\le\frac 12  L_k(t)$,  
\begin{equation}
 \label{reconstr}
 \Big| \lim_{\gl\to 0} \frac{\gs(t)}{g(0)}\int_{|\xi-\phi^\gp_k(t)|<\wt \vep_3} R_{x, \wt \vep_1, \wt \vep_2}^{adp, 2nd,\gl}(t, \xi)d\xi -x_k(t) \Big|\le \Bd_2, 
\end{equation}
where $\Bd_2=\Bd_2^\gp+\Bd_2^{\gp\gp}$ with 
\begin{equation}
\label{def_Bd2}
\begin{array}{l}
\Bd_2^\gp =  \frac{1}{|g(0)|}\big\{ 2\ga_k\big(\wt \vep_1 +\gs(t)\Pi_0(t) 
\big) +{A_k(t)}\; \big|\int_{|u|\ge \ga_k} G_k(u) du\big|+\sum_{\ell\not=k}A_\ell(t) M_{\ell, k}(t)\big\}, \\
 \Bd_2^{\gp \gp} =  \frac{1}{|g(0)|}\big\{ 2\gs(t)\Pi_0(t) 
\big) +\gs(t){A_k(t)}\;\|g\|_1 |Z_t|+\sum_{\ell\not=k}A_\ell(t) M_{\ell, k}(t)\big\}
\end{array}
\end{equation}
and $|Z_t|$ denoting the Lebesgue measure of the set $Z_t$:
 \begin{equation}
 \label{def_Zt}
 Z_t= \{\eta: \; (t, \eta)\in O_k, |V_x(t, \eta)|>\wt \vep_1, \big| \partial _\eta\big({\wt V^{g_1}_x(t, \eta)}/{\wt V_x(t, \eta)}\big)\big| \le \wt \vep_2\}. 
 \end{equation}
\end{theo}

\bigskip

We postpone the proof of Theorem \ref{theo:main} to the end of this section. 


\begin{mrem} 
With the decay conditions of $G_k(u)$ and $G_{j, \ell}(t, \eta)$,  Theorem \ref{theo:main} can be stated in the formulation in Corollary \ref{cor:theo_standand_form}. Here instead of giving such a statement for the 2nd-order adaptive FSST as in Corollary \ref{cor:theo_standand_form}, we look at the estimate bounds when 
$g(t)$ is the Gaussian function given by \eqref{def_g}. 
\end{mrem}

\bigskip

First we look at the bounds for $\Err_1, \Err_2$. From \eqref{def_Err1} and {\eqref{def_Err2}, we have  
\begin{eqnarray*}
|\Err_1|\hskip -0.6 cm && \le 2\pi |B_k(t, \eta)|+ 2\pi \gs(t) |D_k(t, \eta)|
+2\pi \ga_k \Pi_0(t)+  \Pi_2(t) + 2\pi |\phi_k^{\gp\gp}(t)| \gs^2(t)\Pi_1(t),\\
|\Err_2| \hskip -0.6 cm &&\le 4\pi^2\gs(t) |E_k(t, \eta)|+4\pi^2 \gs^2(t) |F_k(t, \eta)|
\\
\nonumber && \qquad 
+ 2\pi \gs(t)\Pi_0(t)+4\pi^2 \ga_k \gs(t) \Pi_1(t) +  2\pi \gs(t) \wt \Pi _1(t)+4\pi^2 
|\phi_k^{\gp\gp}(t)| \gs^3(t)\; \Pi_2(t).
 \end{eqnarray*}
We need to look at the estimates for  $B_k(t, \eta), D_k(t, \eta), E_k(t, \eta), F_k(t, \eta)$, which are 
determined by $G_{j, \ell}(t, \eta)$, for $(t, \eta)\in O_k$. 

When 
$g(t)$ is given by \eqref{def_g}, then $G_k(u)$ is given by \eqref{STFT_LinearChip}, and $\ga_k$ by solving  $|G_k(u)|=\tau_0$ for $u$ is given by \eqref{def_gak_Gaussian}. 
For  $(t, \eta)\in O_k$, we have 
$$
\eta-\phi_{k-1}^\gp(t)\ge \phi_{k}^\gp(t)-\frac{\ga_k}{\gs(t)} -\phi_{k-1}^\gp(t)\ge
 \big(\frac{\ga_{k-1}}{\gs(t)}+\frac{\ga_k}{\gs(t)} \big)- \frac{\ga_k}{\gs(t)}=\frac{\ga_{k-1}}{\gs(t)}. 
$$
Thus 
$$
|G_{0, k-1}(t, \eta)|=|G_{k-1}(\gs(\eta- \phi_{k-1}^\gp(t))|\le |G_{k-1}(\ga_{k-1})|=\tau_0. 
$$
Similarly, we can obtain for  $(t, \eta)\in O_k$
$$
|G_{0, k+1}(t, \eta)| \le \tau_0. 
$$
For other $\ell$, $|G_{0, \ell}(t, \eta)|$ is much smaller than $\tau_0$. 

$G_{j, \ell}(t, \eta)$ can be estimated similarly. For example, for $G_{1, k-1}(t, \eta)$, with \eqref{GjkwithGk}
and 
$$
G_{k-1}^\gp(u)=G_{k-1}(u) (-2u) 
\frac{2\pi^2\big(1+i2\pi \phi_{k-1}^{\gp\gp}(t)\gs^2(t)\big)}{1+\big(2\pi \phi_{k-1}^{\gp\gp}(t)\gs^2(t)\big)^2}, 
$$
we have 
\begin{eqnarray*}
&&|G_{1, k-1}(t, \eta)|=\frac 1{2\pi}|G^\gp_{k-1}\big(\gs(t)(\eta-\phi_{k-1}^\gp(t))\big) |\\
&&=|G_{k-1}\big(\gs(t)(\eta-\phi_{k-1}^\gp(t))\big) |
\frac{ 2\pi |\gs(t) (\eta-\phi_{k-1}^\gp(t))|}{
\sqrt{1+\big(2\pi \phi_{k-1}^{\gp\gp}(t)\gs^2(t)\big)^2}}\\
&&\le |G_{k-1}\big(\ga_{k-1}\big) |  \; 2\pi (\phi_k^\gp(t)-\phi_{k-1}^\gp(t)+\ga_k)=2\pi \tau_0 (\phi_k^\gp(t)-\phi_{k-1}^\gp(t)+\ga_k). 
\end{eqnarray*}
since for $(t, \eta)\in O_k$, $\ga_{k-1}\le \gs(t)(\eta -\phi_{k-1}^\gp(t))
\le \phi_k^\gp(t)-\phi_{k-1}^\gp(t)+\ga_k$. We also can get 
\begin{eqnarray*}
&&|G_{1, k+1}(t, \eta)|\le 2\pi \tau_0(\phi_{k+1}^\gp(t)-\phi_k^\gp(t)+\ga_k). 
\end{eqnarray*}
For other $\ell$, $|G_{1, \ell}(t, \eta)|$ are much smaller than 
$2\pi \tau_0 (\phi_k^\gp(t)-\phi_{k-1}^\gp(t)+\ga_k)$ or $2\pi \tau_0(\phi_{k+1}^\gp(t)-\phi_k^\gp(t)+\ga_k)$. 

For other $j$, one can show that $|G_{j, k\pm1}(t, \eta)|$ are bounded by $C_k \tau_0$, where $C_k$ is a polynomial of $(|\phi_{k\pm1}^\gp(t)-\phi_k^\gp(t)|+\ga_k)$, and $|G_{j, \ell}(t, \eta)|$ for other $\ell $ with $|\ell-k|\ge 2$ are more smaller. 

By the above discussion, we can conclude that for $(t, \eta)\in O_k$, 
$B_k(t, \eta), D_k(t, \eta), E_k(t, \eta), F_k(t, \eta)$ are dominated by $C^\gp_k \tau_0$, where $C^\gp_k$ 
is a polynomial of  $ (\phi_k^\gp(t)-\phi_{k-1}^\gp(t)+\ga_k)$ and $(\phi_{k+1}^\gp(t)-\phi_k^\gp(t)+\ga_k)$
with degree $\le 2$. Therefore, if $\vep_1, \vep_3, \tau_0$ are small, then  $\Bd_1$ in \eqref{def_Bd1} is small.

\bigskip 
Next we look at $\Bd_2$ in \eqref{def_Bd2}. 
First we consider $M_{\ell, k}(t)$ defined by \eqref{def_Mlk}. For $\ell=k-1$, since $\gs(t)(\phi_k^\gp(t)-\phi_{k-1}^\gp(t))> \ga_{k-1}+\ga_k$ and $|G_{k-1}(u)|$ is a decreasing function of $u$, we have 
\begin{eqnarray*}
M_{k-1, k}(t) \hskip -0.6cm &&=  \int_{|u|< \ga_k } | G_{k-1} \big(u+\gs(t)(\phi_k^\gp(t)-\phi_{k-1}^\gp(t) )\big)| du\\
&& \le \int_{|u|< \ga_k } | G_{k-1} (u+\ga_{k-1}+\ga_k )| du \\
&& \le  2 \ga_k  |G_{k-1}(\ga_{k-1})| =2\ga_k \tau_0.  
\end{eqnarray*}
We can show similarly  $M_{k+1, k}(t)\le 2\ga_k \tau_0$. For other $\ell$, $M_{\ell, k}(t)$ will be much smaller than $2\ga_k \tau_0$. 

\bigskip 
Now let us look at $\int_{|u|\ge \ga_k} G_k(u) du$: 
\begin{eqnarray*}
&&|\int_{|u|\ge \ga_k} G_k(u) du| \le \int_{|u|\ge \ga_k} |G_k(u)| du=\frac 1{\big(1+(2\pi \phi_k^{\gp\gp}(t)\gs^2(t))^2\big)^{\frac 14}}\;
 \int_{|u|\ge \ga_k} e^{-\frac{2\pi^2}{1+(2\pi \phi_k^{\gp\gp}(t)\gs^2(t))^2}u^2}du\\
 &&=\big(1+(2\pi \phi_k^{\gp\gp}(t)\gs^2(t))^2\big)^{\frac 14}\;
 \int_{|u|\ge \frac{\ga_k}{\sqrt{1+(2\pi \phi_k^{\gp\gp}(t)\gs^2(t))^2}} 
 } e^{-2\pi^2 u^2}du\\
 &&\le 
 \big(1+(2\pi \phi_k^{\gp\gp}(t)\gs^2(t))^2\big)^{\frac 14}\;
 \int_{|u|\ge \ga
 } e^{-2\pi^2 u^2}du \le   \big(1+(2\pi \phi_k^{\gp\gp}(t)\gs^2(t))^2\big)^{\frac 14}\; 
 \frac 1{\sqrt {2\pi}}\frac{\tau_0}{1+\sqrt{1-\tau_0}}, 
\end{eqnarray*}
where the last inequality follows from \eqref{est_Guasian_intgral}. 

\bigskip 

Finally we present the proof of Theorem \ref{theo:main}. 

{\bf Proof  of  Theorem \ref{theo:main} Part (a)}.  Assume $(t, \eta)\not \in \cup _{k=1}^K O_k$. Then for any $k$, by the definition of $O_k$ in \eqref{def_Ok}, we have  $|G_k\big(\gs(t)(\eta-\phi_k^\gp(t) )\big)|\le \tau_0$. Hence, by \eqref{STFT_approx} and \eqref{err0_est}, we have 
\begin{eqnarray*}
|\wt V_x(t, \eta)|\hskip -0.6cm &&\le |\err_0|+\sum_{k=1}^K |x_k(t) G_k\big(\gs(t)(\eta-\phi_k^\gp(t) )\big)| 
\\
&&\le  \gs(t) \Pi_0(t)+\tau _0 \sum_{k=1}^K A_k(t)\le \wt \vep_1, 
\end{eqnarray*}
a contradiction to the assumption $|\wt V_x(t, \eta)|>\wt \vep_1$. Hence the statement in (a) holds. 
\hfill $\blacksquare$

\bigskip 

{\bf Proof  of  Theorem \ref{theo:main} Part ${\rm (b)}$}.  Plugging  $\pd_t \wt V_x(t, \eta)$ in \eqref{result_lem1} to $\go_{x}^{adp, 2nd, c}$ in \eqref{def_transformation_2nd_complex}, we have 
\begin{eqnarray*}
&&\go_{x}^{adp, 2nd, c}=\frac {{\partial}_t \wt V_x(t, \eta)}{i2\pi \wt V_x(t, \eta)}+\frac{\gs'(t)}{i2\pi \gs(t)}
- \frac{\wt V^{g_1}_x(t, \eta)}{i2\pi \wt V_x(t, \eta)} P_0(t, \eta)
+ \frac {\gs'(t)}{\gs(t)}\frac {\wt V^{g_3}_x(t, \eta)}{i2\pi \wt V_x(t, \eta)}\\
&&=\frac 1{i2\pi \wt V_x(t, \eta)}\Big\{\big(i2\pi \phi^\gp_k(t)-\frac {\gs'(t)}{\gs(t)}\big)\wt V_x(t, \eta)+
i2\pi \phi_k^{\prime\prime}(t)\gs(t)\wt V^{g_1}_x(t, \eta)- \frac {\gs'(t)}{\gs(t)}\wt V^{g_3}_x(t, \eta)+\Err_1\Big\}
\\
&&\qquad + \frac{\gs'(t)}{i2\pi \gs(t)} - \frac{\wt V^{g_1}_x(t, \eta)}{i2\pi \wt V_x(t, \eta)} P_0(t, \eta)
+ \frac {\gs'(t)}{\gs(t)}\frac {\wt V^{g_3}_x(t, \eta)}{i2\pi \wt V_x(t, \eta)} 
\\
&&=\phi^\gp_k(t)+ \phi_k^{\prime\prime}(t)\gs(t)\frac{\wt V^{g_1}_x(t, \eta)}{\wt V_x(t, \eta)}+\frac{\Err_1}{i2\pi \wt V_x(t, \eta)} - \frac{\wt V^{g_1}_x(t, \eta)}{i2\pi\wt V_x(t, \eta)} P_0(t, \eta)\\
&&=\phi^\gp_k(t)+ \phi_k^{\prime\prime}(t)\gs(t)\frac{\wt V^{g_1}_x(t, \eta)}{\wt V_x(t, \eta)}+\frac{\Err_1}{i2\pi \wt V_x(t, \eta)} - \frac{\wt V^{g_1}_x(t, \eta)}{i2\pi\wt V_x(t, \eta)} \big(i2\pi \gs(t)\phi^{\prime\prime}_k(t)+\Err_3\big)\\
&&=\phi^\gp_k(t)+\frac{\Err_1}{i2\pi \wt V_x(t, \eta)}- \frac{\wt V^{g_1}_x(t, \eta)\Err_3}{i2\pi\wt V_x(t, \eta)} \\
&&=\phi^\gp_k(t)+\Err_4,  
\end{eqnarray*}
where the last third equation follows from \eqref{result_lem3}. 
This shows \eqref{transformation_approx_2nd}.

For \eqref{IF_est_Bd1}, with the assumptions   
$|\wt V_x(t, \eta)|>\wt \vep_1$ and 
$$
|\pd _\eta \big(\wt V^{g_1}_x(t, \eta)/\wt V_x(t, \eta)\big)|
=\Big|\wt V_x(t, \eta)\pd_\eta \wt V^{g_1}_x(t, \eta)-\pd_\eta \wt V_x(t, \eta) \wt V^{g_1}_x(t, \eta)\Big|\big/|\wt V_x(t, \eta)|^2
>\wt \vep_2, 
$$ 
we have 
\begin{eqnarray}
\nonumber 
&&|\Err_4|=\Big| \frac{\Err_1}{i2\pi \wt V_x(t, \eta)}- \frac{\wt V^{g_1}_x(t, \eta)\Err_3}{i2\pi\wt V_x(t, \eta)}
\Big|\\
\nonumber &&=\Big| \frac{\Err_1}{i2\pi \wt V_x(t, \eta)}- \frac{\wt V^{g_1}_x(t, \eta)}{i2\pi\wt V_x(t, \eta)}
\frac {\wt V_x(t, \eta) \Err_2-\pd_\eta \wt V_x(t, \eta) \Err_1 }{\wt V_x(t, \eta)\pd_\eta 
\wt V^{g_1}_x(t, \eta)-\wt V^{g_1}_x(t, \eta)\pd_\eta 
\wt V_x(t, \eta)}
\Big|\\
\nonumber &&\le \frac{|\Err_1|}{2\pi |\wt V_x(t, \eta)|}+ \frac{|\wt V^{g_1}_x(t, \eta)|}{2\pi|\wt V_x(t, \eta)|}
\frac {\Big(|\pd_\eta \wt V_x(t, \eta)| \; |\Err_1| + |\wt V_x(t, \eta)| \; |\Err_2|\Big)\big /|\wt V_x(t, \eta)|^2}{
\Big|\wt V_x(t, \eta)\pd_\eta 
\wt V^{g_1}_x(t, \eta)-\wt V^{g_1}_x(t, \eta)\pd_\eta 
\wt V_x(t, \eta)\Big| \big/|\wt V_x(t, \eta)|^2}\\
&&< \frac{|\Err_1|}{2\pi \wt \vep_1}+ \frac1{2\pi \wt \vep_1^3 \wt \vep_2}
|\wt V^{g_1}_x(t, \eta)| \big(|\pd_\eta \wt V_x(t, \eta)|\;  |\Err_1| + \wt \vep_1 \; |\Err_2|\big). 
\label{est_Res4} 
\end{eqnarray}
Thus $|\Err_4|< \Bd_1$, as desired. \hfill $\blacksquare$

\bigskip 

{\bf Proof  of  Theorem \ref{theo:main} Part (c)}. First we have the following result which can be derived as that on p.254 in \cite{Daub_Lu_Wu11}:
\begin{equation}
\label{FSST_STFT_relation1}
 \lim_{\gl\to 0}\int_{|\xi-\phi^\gp_k(t)|<\wt \vep_3} R_{x, \wt \vep_1, \wt \vep_2}^{adp, 2nd,\gl}(t, \xi)d\xi =
 \int_{X_t}  \wt V_x(t, \eta) 
 d\eta ,
\end{equation}
where 
\begin{equation*}
X_t=\big\{\eta:\; |\wt V_x(t, \eta)|>\wt \vep_1,  \;  \big| \partial _\eta\big({\wt V^{g_1}_x(t, \eta)}/{\wt V_x(t, \eta)}\big)\big| > \wt \vep_2 \; 
\hbox{and $\big|\phi^\gp_k(t)-\go_{x, \wt \vep_1, \wt \vep_2}^{adp, 2nd}(t, \eta)\big|<\wt \vep_3$}\big \}. 
\end{equation*}

Denote 
\begin{equation*}
Y_t=\big\{\eta:\; |\wt V_x(t, \eta)|>\wt \vep_1,  \; \big| \partial _\eta\big({\wt V^{g_1}_x(t, \eta)}/{\wt V_x(t, \eta)}\big)\big| > \wt \vep_2 \; \hbox{and $(t, \eta)\in O_k$}\big \}. 
\end{equation*}
Then we have $Y_t=X_t$. Indeed, 
by Theorem \ref{theo:main} part (b), if $\eta\in Y_t$, then $\big|\phi^\gp_k(t)-\go_{x, \wt \vep_1, \wt \vep_2}^{adp, 2nd}(t, \eta)\big|<\Bd_1\le \wt \vep_3$. Thus $\eta\in X_t$. Hence $Y_t\subseteq X_t$. 

On the other hand, suppose $\eta\in X_t$. Since $|\wt V_x(t, \eta)|>\wt \vep_1$, by Theorem \ref{theo:main} part (a), $(t, \eta)\in O_\ell$ for an $\ell$ in $\{1, 2, \cdots, K\}$. If $\ell\not =k$, then 
\begin{eqnarray*}
\big|\phi^\gp_k(t)-\go_{x, \wt \vep_1, \wt \vep_2}^{adp, 2nd}(t, \eta)\big|\hskip -0.6cm &&\ge |\phi^\gp_k(t) -\phi^\gp_\ell(t)|-\big|\phi^\gp_\ell(t)-\go_{x, \wt \vep_1, \wt \vep_2}^{adp, 2nd}(t, \eta)\big|\\
&&> L_k(t)-\wt \vep_3\ge \wt \vep_3,  
\end{eqnarray*}
and this contradictes to the assumption $\eta\in X_t$ with $\big|\phi^\gp_k(t)-\go_{x, \wt \vep_1, \wt \vep_2}^{adp, 2nd}(t, \eta)\big|<\wt \vep_3$, where we have used the fact $|\phi^\gp_k(t) -\phi^\gp_\ell(t)|\ge L_k(t)$ and 
$\big|\phi^\gp_\ell(t)-\go_{x, \wt \vep_1, \wt \vep_2}^{adp, 2nd}(t, \eta)\big|<\Bd_1\le \wt \vep_3$ by Theorem \ref{theo:main} part (b). 
Hence $\ell=k$ and $\eta\in Y_t$. Thus we know $X_t=Y_t$. 

The facts $X_t=Y_t$ and $Y_t\cap Z_t=\emptyset$, togehter with \eqref{FSST_STFT_relation1}, imply that  
\begin{eqnarray}
 \nonumber &&\lim_{\gl\to 0}\int_{|\xi-\phi^\gp_k(t)|<\wt \vep_3} R_{x, \wt \vep_1, \wt \vep_2}^{adp, 2nd,\gl}(t, \xi)d\xi =
 \int_{Y_t}\wt V_x(t, \eta) d\eta  =\int_{Y_t\cup Z_t}\wt V_x(t, \eta) d\eta - \int_{Z_t}\wt V_x(t, \eta) d\eta \\
\label{FSST_STFT_relation2} && = \int_{\{|\wt V_x(t, \eta)|>\wt \vep_1\}\cap \{\eta: (t, \eta)\in O_k\}}  \wt V_x(t, \eta) 
 d\eta - \int_{Z_t}\wt V_x(t, \eta) d\eta. 
\end{eqnarray}

Furthermore, 
\begin{eqnarray*}
&&\big|\gs(t) \int_{\{|\wt V_x(t, \eta)|>\wt \vep_1\}\cap \{\eta: (t, \eta)\in O_k\}}  \wt V_x(t, \eta) d\eta 
-g(0) x_k(t)\big|\\
&&= \big| \gs(t)\int_{\{\eta: (t, \eta)\in O_k\}}  \wt V_x(t, \eta) d\eta 
-g(0) x_k(t)- \gs(t)\int_{\{|\wt V_x(t, \eta)|\le\wt \vep_1\}\cap \{\eta: (t, \eta)\in O_k\}}  \wt V_x(t, \eta) d\eta 
\big|\\
&&\le  \gs(t)\wt \vep_1 \; \frac{2\ga_k}{\gs(t)}+\big|\gs(t) \int_{\{\eta: (t, \eta)\in O_k\}}  \big(\sum_{\ell=1}^K x_\ell(t) G_{0, \ell}(t, \eta)+\err_0 \big)d\eta  
-g(0) x_k(t)\big| \\
&& \le   2\wt \vep_1 \ga_k+
 \gs(t) |\err_0| \; \frac{2\ga_k}{\gs(t)}+
\big| x_k(t) \int_{|u|<\ga_k}   G_{k}(u)du   
-g(0) x_k(t)\big| \\
&&\qquad +\sum_{\ell\not=k }A_\ell(t) |\gs(t)\int_{\{\eta: (t, \eta)\in O_k\}} G_{0, \ell}(t, \eta) d\eta| 
\\
&&\le 2(|\err_0|+\wt \vep_1) \ga_k+
\big|x_k(t) \int_\R   G_{k}(u)du  
-g(0) x_k(t)-x_k(t)\int_{|u|\ge \ga_k}   G_{k}(u)du  \big|+\sum_{\ell\not=k }A_\ell(t) M_{\ell, k}(t)\\
&&= 2(|\err_0|+\wt \vep_1)\ga_k  +
\big|x_k(t) g(0)-g(0) x_k(t)-x_k(t)\int_{|u|\ge \ga_k}   G_{k}(u)du
 \big|+\sum_{\ell\not=k }A_\ell(t) M_{\ell, k}(t) \\
&&=  2(|\err_0|+\wt \vep_1)\ga_k+A_k(t)\; \big|\int_{|u|\ge \ga_k}   G_{k}(u)du \big|+\sum_{\ell\not=k }A_\ell(t) M_{\ell, k}(t), 
\end{eqnarray*}
where we have used the fact:
\begin{eqnarray*}
&&
\int_\R G_{k}(u)du =\int_\R {\cal F}\big\{e^{i\pi \phi_k^{\gp\gp}(t)\tau^2}g(\tau)\big\}(u) du
=\big(e^{i\pi \phi_k^{\gp\gp}(t)\tau^2}g(\tau)\big)\Big|_{\tau=0}=g(0). 
\end{eqnarray*}
Hence, we have 
\begin{equation}
\label{est_Bd3p}
\big|\frac {\gs(t)}{g(0)} \int_{\{|\wt V_x(t, \eta)|>\wt \vep_1\}\cap \{\eta: (t, \eta)\in O_k\}}  \wt V_x(t, \eta) d\eta 
-x_k(t)\big|\le \Bd_2^\gp. 
\end{equation}

In addition, 
\begin{eqnarray*}
&&\big|\gs(t) \int_{Z_t}  \wt V_x(t, \eta) d\eta\big|=
\gs(t) \big| \int_{Z_t}  \big(\sum_{\ell=1}^K x_\ell(t) G_{0, \ell}(t, \eta)+\err_0 \big)d\eta \big| \\
&& \le    \gs(t) |\err_0| \; \frac{2\ga_k}{\gs(t)}+\gs(t)A_k(t) \sup_{\eta \in Z_t} 
|G_k\big(\gs(t)(\eta-\phi_k^{\rq{}}(t)\big)| \; |Z_t|+
\sum_{\ell\not=k }A_\ell(t) |\gs(t)\int_{\{\eta: (t, \eta)\in O_k\}} G_{0, \ell}(t, \eta) d\eta| 
\\
&&\le 2 |\err_0| \ga_k+\gs(t)A_k(t) \|g\|_1 \; |Z_t|+\sum_{\ell\not=k }A_\ell(t) M_{\ell, k}(t)
\le \Bd_2^{\gp\gp}. 
\end{eqnarray*}
The above estimates, together with \eqref{FSST_STFT_relation2}, leads to \eqref{reconstr}. 
This completes the proof of Theorem \ref{theo:main} part (c). 
\hfill $\blacksquare$

\section*{Appendix}
In this appendix, we provide the proof of Lemmas \ref{lem:lem1} and \ref{lem:lem3}. For simplicity of presentation, we drop $x, t, \eta$ in  $\wt V_x(t, \eta)$. 

{\bf Proof  of  Lemma \ref{lem:lem1}}. \quad 
By \eqref{STFT_approx0} with $g$ replaced by $g^\gp$, 
\begin{eqnarray*}
&&\wt V^{g^\gp}=\sum_{\ell=1}^K \int_\R x_\ell(t)e^{i2\pi (\phi_\ell^\gp(t) \tau +\frac12\phi^{\gp\gp}_\ell(t) \tau^2)}\frac 1{\gs(t)}g^\gp(\frac \tau{\gs(t)}) e^{-i2\pi \eta \tau}d\tau +\err_0^\gp\\
&&=\sum_{\ell=1}^K \int_\R x_\ell(t)e^{-i2\pi (\eta-\phi_\ell^\gp(t)) \tau +i\pi \phi^{\gp\gp}_\ell(t) \tau^2} \frac {\pd}{\pd\tau} \Big(g(\frac \tau{\gs(t)})\Big) d\tau +\err_0^\gp\\
&&=-\sum_{\ell=1}^K \int_\R \frac {\pd}{\pd\tau} \Big( x_\ell(t)e^{-i2\pi (\eta-\phi_\ell^\gp (t) )\tau +i\pi \phi^{\gp\gp}_\ell(t) \tau^2}\Big)g(\frac \tau{\gs(t)}) d\tau +\err_0^\gp\\
&&=i2\pi \sum_{\ell=1}^K x_\ell(t)(\eta-\phi_\ell^\gp (t))\int_\R  e^{-i2\pi (\eta-\phi_\ell^\gp (t) )\tau +i\pi \phi^{\gp\gp}_\ell(t) \tau^2}g(\frac \tau{\gs(t)}) d\tau\\
&&\quad  -
i2\pi \sum_{\ell=1}^K x_\ell(t)\phi_\ell^{\gp\gp} (t) 
\int_\R  e^{-i2\pi (\eta-\phi_\ell^\gp (t))\tau +i\pi \phi^{\gp\gp}_\ell(t) \tau^2}\tau g(\frac \tau{\gs(t)}) d\tau+\err_0^\gp\\
&&=i2\pi \gs(t)\sum_{\ell=1}^K x_\ell(t)(\eta-\phi_\ell^\gp (t))G_{0, \ell}(t, \eta) 
 - i2\pi \gs^2(t)\sum_{\ell=1}^K x_\ell(t)\phi_\ell^{\gp\gp} (t)G_{1, \ell}(t, \eta)  +\err_0^\gp. 
\end{eqnarray*}
This and \eqref{result_partial_tV} imply that 
\begin{eqnarray*}
&&\pd_t \wt V+\frac {\gs'(t)}{\gs(t)}(\wt V+\wt V^{g_3})
-i2\pi \phi^\gp_k(t)\wt V-i2\pi \phi_k^{\prime\prime}(t)\gs(t)\wt V^{g_1}\\
&&=i2\pi \eta \wt V-\frac 1{\gs(t)}\wt V^{g^\gp}-i2\pi \phi^\gp_k(t)\wt V-i2\pi \phi_k^{\prime\prime}(t)\gs(t)\wt V^{g_1}\\
&&=i2\pi \big(\eta -\phi_k^\gp(t)\big)\Big (\sum_{\ell=1}^K x_\ell(t)G_{0, \ell}(t, \eta)+\err_0\Big)\\
&&\qquad -\frac 1{\gs(t)}\Big(i2\pi \gs(t)\sum_{\ell=1}^K x_\ell(t)(\eta-\phi_\ell^\gp (t))G_{0, \ell}(t, \eta) 
 - i2\pi \gs^2(t)\sum_{\ell=1}^K x_\ell(t)\phi_\ell^{\gp\gp} (t)G_{1, \ell}(t, \eta)  +\err_0^\gp\Big)\\
 &&\qquad 
-i2\pi \phi_k^{\prime\prime}(t)\gs(t)\Big (\sum_{\ell=1}^K x_\ell(t)G_{1, \ell}(t, \eta)+\err_1\Big)
\\
&&=i2\pi \sum_{\ell\not=k} x_\ell(t)(\phi_\ell^\gp (t)-\phi_k^\gp (t))G_{0, \ell}(t, \eta) +
i2\pi \gs(t)\sum_{\ell\not= k} x_\ell(t)(\phi_\ell^{\gp\gp} (t)-\phi_k^{\gp\gp} (t))G_{1, \ell}(t, \eta)\\
&&\qquad +i2\pi \big(\eta -\phi_k^\gp(t)\big)\err_0-\frac{\err_0^\gp}{\gs(t)}-i2\pi \phi_k^{\prime\prime}(t) \gs(t) \; \err_1\\
&&=i2\pi B_k(t, \eta)+i2\pi \gs(t) D_k(t, \eta)+i2\pi \big(\eta -\phi_k^\gp(t)\big)\err_0-\frac{\err_0^\gp}{\gs(t)}-i2\pi \phi_k^{\prime\prime}(t)\gs(t)\;  \err_1\\
&&=\Err_1.  
\end{eqnarray*}
This completes the proof of Lemma \ref{lem:lem1}. \hfill $\blacksquare$ 

\bigskip 

{\bf Proof of Lemma \ref{lem:lem3}}. \quad 
 One can obtain for $j\ge 1$, 
\begin{eqnarray*}
\frac{\partial}{\partial \eta} G_{j, k}(t, \eta)=-i2\pi \gs(t) G_{j+1, k}(t, \eta). 
\end{eqnarray*}
Thus we have 
$$
\partial_\eta B_k(t, \eta)=-i2\pi \gs(t) E_k(t, \eta), \; \partial_\eta D_k(t, \eta)=-i2\pi \gs(t) F_k(t, \eta). 
$$
In addition, it is straightforward to verify that 
$$
\partial_\eta  \err_j=-i2\pi \gs(t)\;  \err_{j+1}, \; \partial_\eta  \err^\gp_j=-i2\pi \gs(t)\;  \err^\gp_{j+1}. 
$$
Hence we have 
\begin{equation}
\label{Err1_Err2}
\partial_\eta  \Err_1= \Err_2. 
\end{equation}

Taking the partial derivate with respect to $\eta$ to the both sides of  \eqref{result_lem1}  and using \eqref{Err1_Err2}, we have 
\begin{equation}
\label{result_lem2}
\pd_\eta \pd_t \wt V=\big(i2\pi \phi^\gp_k(t)-\frac {\gs'(t)}{\gs(t)}\big)\pd_\eta 
\wt V+i2\pi \phi_k^{\prime\prime}(t)\gs(t)\pd_\eta \wt V^{g_1}
- \frac {\gs'(t)}{\gs(t)}\pd_\eta \wt V^{g_3}+\Err_2. 
\end{equation}
Note that 
$$
P_0(t, \eta)=\frac1{\wt V\pd_\eta 
\wt V^{g_1}-\wt V^{g_1}\pd_\eta 
\wt V}\Big(\wt V\pd_\eta\pd_t 
\wt V-\pd_\eta\wt V\pd_t 
\wt V+\frac {\gs^\gp(t)}{\gs(t)}(\wt V\pd_\eta 
\wt V^{g_3}-\wt V^{g_3}\pd_\eta 
\wt V) \Big). 
$$
Thus, by \eqref{result_lem1} and \eqref{result_lem2}, 
\begin{eqnarray*}
&&\big(P_0(t, \eta)-i2\pi \gs(t)\phi^{\prime\prime}_k(t)\big)\big(\wt V\pd_\eta 
\wt V^{g_1}-\wt V^{g_1}\pd_\eta 
\wt V)\\
&&=\wt V\pd_\eta\pd_t 
\wt V-\pd_\eta\wt V^{g_1}\pd_t 
\wt V+\frac {\gs^\gp(t)}{\gs(t)}(\wt V\pd_\eta 
\wt V^{g_3}-\wt V^{g_3}\pd_\eta 
\wt V) -i2\pi \gs(t)\phi^{\prime\prime}_k(t)\big(\wt V\pd_\eta 
\wt V^{g_1}-\wt V^{g_1}\pd_\eta 
\wt V)\\
&&=\wt V\Big(
\big(i2\pi \phi^\gp_k(t)-\frac {\gs'(t)}{\gs(t)}\big)\pd_\eta 
\wt V+i2\pi \phi_k^{\prime\prime}(t)\gs(t)\pd_\eta \wt V^{g_1}
- \frac {\gs'(t)}{\gs(t)}\pd_\eta \wt V^{g_3}+\Err_2
\Big)\\
&&\qquad -\pd_\eta\wt V\Big(
\big(i2\pi \phi^\gp_k(t)-\frac {\gs'(t)}{\gs(t)}\big)\wt V+
i2\pi \phi_k^{\prime\prime}(t)\gs(t)\wt V^{g_1}- \frac {\gs'(t)}{\gs(t)}\wt V^{g_3}+\Err_1
\Big)\\
&&\qquad +\frac {\gs^\gp(t)}{\gs(t)}(\wt V\pd_\eta 
\wt V^{g_3}-\wt V^{g_3}\pd_\eta 
\wt V) -i2\pi \gs(t)\phi^{\prime\prime}_k(t)\big(\wt V\pd_\eta 
\wt V^{g_1}-\wt V^{g_1}\pd_\eta 
\wt V)\\
&&=\wt V\; \Err_2 -\pd_\eta\wt V\; \Err_1
\end{eqnarray*}
Therefore,  we have 
$$
P_0(t, \eta)-i2\pi \gs(t)\phi^{\prime\prime}_k(t)=
\frac{\wt V\; \Err_2 -\pd_\eta\wt V\; \Err_1}{\wt V\pd_\eta 
\wt V^{g_1}-\wt V^{g_1}\pd_\eta 
\wt V}=\Err_3, 
$$
as desired.  This completes the proof of Lemma \ref{lem:lem3}. \hfill $\blacksquare$ 

{\bf Acknowledgments:}  Any opinions, findings and conclusions or recommendations expressed in this material are those of the authors and do not necessarily reflect the views of AFRL (Air Force Research Laboratory). 
Q. Jiang was supported by  the 2018 Air Force Visiting Faculty Research Program (VFRP) and by Simons Foundation (Grant No. 353185).


\begin{thebibliography}{10}

\bibitem{A_Flandrin_reassignment95}  F. Auger and P. Flandrin, Improving the readability of time-frequency and  time-frequency representations by the reassignment method,  {IEEE Trans. Signal Proc.}, 43 (1995), pp. 1068--1089.

\bibitem{BMO18} R. Behera, S. Meignen, and T. Oberlin, Theoretical analysis of the 2nd-order synchrosqueezing transform, {Appl. Comput. Harmon. Anal.}, 45 (2018), pp. 374--404. 


\bibitem{Berrian18} A.J.  Berrian,  {The Chirped Quilted Synchrosqueezing Transform and its Application to Bioacoustic Signal Analysis}, Ph.D. Dissertation, University of California at Davis,  2018.

\bibitem{Saito17} A.J.  Berrian and N. Saito, Adaptive synchrosqueezing based on a quilted short-time Fourier transform,
arXiv:1707.03138v5, Sep. 2017. 


\bibitem{Chui_Lin_Wu15} C.K. Chui, Y.-T. Lin, and H.-T. Wu,  Real-time dynamics acquisition from irregular samples - with application to anesthesia evaluation, {Anal. Appl.}, 14 (2016), pp.537--590.


\bibitem{Chui_Mhaskar15}  C.K. Chui and H.N. Mhaskar, Signal decomposition and analysis via extraction of frequencies, {Appl. Comput. Harmon. Anal.}, 40 (2016), pp. 97--136.

\bibitem{Chui_Walt15} C.K. Chui and M.D. van der Walt, Signal analysis via instantaneous frequency estimation of signal components,  {Int'l  J.  Geomath.}, 6 (2015), pp. 1--42.


\bibitem{Daub_Lu_Wu11} I. Daubechies, J. Lu, and H.-T. Wu, Synchrosqueezed wavelet transforms:
An empirical mode decomposition-like tool,  {Appl. Comput. Harmon. Anal.}, 30 (2011), pp. 243--261.

\bibitem{Daub_Wang_Wu15} I. Daubechies, Y. Wang, and H.-T. Wu, ConceFT: Concentration of frequency and time via a multitapered synchrosqueezed transform,
{Phil. Trans. Royal Soc. A,} 374 (2016). 

\bibitem{Huang98}  N.E. Huang, Z. Shen, S.R. Long, M.L. Wu, H.H. Shih, Q. Zheng, N.C. Yen, C.C. Tung,  and H.H. Liu, The empirical mode decomposition and Hilbert spectrum for nonlinear and nonstationary
time series analysis,  {Proc. Roy. Soc. London A}, 454 (1998), pp. 903--995.

 \bibitem{Iatsenko15} D. Iatsenko, P.-V. E. McClintock, and A. Stefanovska, Linear and synchrosqueezed time-frequency representations revisited: Overview, standards of use, resolution, reconstruction, concentration, and algorithms, {Digital Signal Proc.},  42 (2015),  pp. 1--26.


\bibitem{Jiang_Suter17} Q.T. Jiang and B.W.  Suter,  Instantaneous frequency estimation based on synchrosqueezing wavelet transform, {Signal Proc.}, 138 (2017),  167--181.


\bibitem{LCHJJ18}  L. Li, H.Y. Cai, H.X. Han, Q.T. Jiang and H.B. Ji, Adaptive short-time Fourier transform and synchrosqueezing transform for non-stationary signal separation, preprint,  2018.

\bibitem{LCJJ17}  L. Li, H.Y. Cai, Q.T.  Jiang and H.B.  Ji,  An empirical signal separation algorithm based on linear time-frequency analysis, Mechanical Systems and Signal Proc., 121 (2019), 791--809. 

\bibitem{LCJJ18}  L. Li, H.Y. Cai, Q.T. Jiang and H.B. Ji, Adaptive synchrosqueezing transform with a time-varying parameter for non-stationary signal separation, preprint, 2018.

\bibitem{Li_Liang12} C. Li and M. Liang, A generalized synchrosqueezing transform for
enhancing signal time-frequency representation,  {Signal Proc.}, 92 (2012), pp. 2264--2274.





















\bibitem{Meignen17} S. Meignen, D.-H. Pham, and S. McLaughlin, On demodulation, ridge detection and synchrosqueezing for multicomponent signals, {IEEE Trans. Signal Proc.}, 65 (2017), pp. 2093--2103.

\bibitem{OM17} T. Oberlin and S. Meignen, The 2nd-order wavelet synchrosqueezing transform,  in 2017 IEEE International Conference on Acoustics, Speech and Signal Processing (ICASSP), March 2017, New Orleans, LA, USA.

\bibitem{MOM14} T. Oberlin, S. Meignen, and V. Perrier, The Fourier-based synchrosqueezing
transform,  in  Proc. 39th Int. Conf. Acoust., Speech,
Signal Proc. (ICASSP), 2014, pp. 315--319.



 \bibitem{MOM15} T. Oberlin, S. Meignen, and V. Perrier, Second-order synchrosqueezing transform or invertible reassignment? Towards ideal time-frequency representations,  {IEEE Trans. Signal Proc.}, 63 (2015),  pp.1335--1344.

\bibitem{Pham17} 	D.-H. Pham and S. Meignen, High-order synchrosqueezing transform for multicomponent signals analysis - With an application to gravitational-wave signal, {IEEE Trans. Signal Proc.}, 65 (2017), pp. 3168--3178.

\bibitem{Wu17} Y.-L. Sheu, L.-Y. Hsu, P.-T. Chou, and H.-T. Wu, Entropy-based time-varying window width selection for nonlinear-type time-frequency analysis, {Int'l J. Data Sci. Anal.}, 3 (2017),  pp. 231--245.


\bibitem{Thakur_Wu11} G. Thakur and H.-T. Wu, Synchrosqueezing based recovery of instantaneous frequency from nonuniform samples, {SIAM J. Math. Anal.}, 43 (2011), pp. 2078--2095.

\bibitem{Thakur_etal_Wu13} G. Thakur, E. Brevdo, N. Fu$\check{\rm c}$kar, and H.-T. Wu, The synchrosqueezing algorithm for time-varying spectral analysis: Robustness properties and new paleoclimate applications,  {Signal Proc.}, 93 (2013), pp. 1079--1094.


\bibitem{Wang_etal14} S.B. Wang, X.F. Chen, G.G. Cai, B.Q. Chen, X. Li, and Z.J. He, Matching demodulation
transform and synchrosqueezing in time-frequency analysis, {IEEE Trans. Signal Proc.}, 62 (2014), pp. 69--84.


\bibitem{WCSGTZ18} S.B. Wang, X.F. Chen, I.W. Selesnick, Y.J. Guo,  C.W. Tong and X.W. Zhang, 
Matching synchrosqueezing transform: A useful tool for characterizing signals with fast varying instantaneous frequency and application to machine fault diagnosis, Mechanical Systems and Signal Proc., 100 (2018),  pp. 242--288.

\bibitem{Wu_thesis} H.-T. Wu, {Adaptive Analysis of Complex Data Sets},  Ph.D. dissertation, Princeton Univ., Princeton, NJ,  2012.

\bibitem{Yang15} H.Z. Yang, Synchrosqueezed wave packet transforms and diffeomorphism based spectral analysis for 1D general mode decompositions, {Appl. Comput. Harmon. Anal.}, 39 (2015), pp.33--66.

\bibitem{Yang18} H.Z. Yang, Statistical analysis of synchrosqueezed transforms, {Appl. Comput. Harmon. Anal.}, 45 (2018), pp.526--550.



























\end{thebibliography}
\end{document}